\documentclass{amsart}

\usepackage{amsfonts, amsbsy, amsmath, amsthm, amssymb, latexsym, enumerate,verbatim,mathrsfs}
\usepackage{color}

\newtheorem{propo}{Proposition}[section]
\newtheorem{conje}{Conjecture}
\newtheorem{lemma}[propo]{Lemma}
\newtheorem{corol}[propo]{Corollary}

\newtheorem{theo}[propo]{Theorem}

\theoremstyle{definition}

\newtheorem{rem}{Remark}

\numberwithin{equation}{section}

\newcommand{\bl}{\begin{lemma}\label}
\newcommand{\el}{\end{lemma}}

\newcommand{\ld}{,\ldots ,}

\newcommand{\ra}{ \rightarrow }
\newcommand{\lan}{ \langle }
\newcommand{\ran}{ \rangle }

\newcommand{\diag}{\mathop{\rm diag}\nolimits}

\newcommand{\Id}{\mathop{\rm Id}\nolimits}

\newcommand{\Jor}{\mathop{\rm Jord}\nolimits}

\newcommand{\FF}{\mathbb{F}}
\newcommand{\GG}{\mathop{\bf G}\nolimits}

\newcommand{\al}{\alpha}
\newcommand{\be}{\beta}

\newcommand{\lam}{\lambda }

\newcommand{\Om}{\Omega }

\newcommand{\up}{^{-1}}
\newcommand{\un}{\mathsf{u}}

\newcommand{\si}{\sigma }
\def\d12{{_{12}}}

\def\acf{{algebraically closed field }}

\def\au{{automorphism }}

\def\f{{following }}

\def\ii{{if and only if }}

\def\ir{{irreducible }}

\def\irt{{irreducible. }}
\def\irr{{irreducible representation }}

\def\itf{{It follows that }}

\def\rep{{representation }}

\def\syl{{Sylow $p$-subgroup }}

\def\SL{{\rm SL}}
\def\SU{{\rm SU}}

\def\Aut{{\rm Aut}}

\def\2k{2^k}
\newcommand{\bc}{\begin{conje}}
\newcommand{\enc}{\end{conje}}

\newcommand{\bp}{\begin{proof} }
\newcommand{\enp}{\end{proof}}

\usepackage[table,dvipsnames]{xcolor}

\usepackage{hyperref}
\hypersetup{ 
 colorlinks =true,
 linkcolor=Sepia,
 citecolor=Plum,
 urlcolor=olive,
 }

\begin{document}
 
\title[Minimal generation of finite simple groups of Lie type by regular\ldots]
{Minimal generation of finite simple groups of Lie type by regular unipotent elements}
 
\author{M.A. Pellegrini}
 
\address{Dipartimento di Matematica e Fisica, Universit\`a Cattolica del Sacro Cuore,
Via della Garzetta 48, 25133 Brescia, Italy}
\email{marcoantonio.pellegrini@unicatt.it}

\author{A.E. Zalesski}
\address{
Universidade de Bras\'ilia, Bras\'ilia-DF, Brazil}
\email{alexandre.zalesski@gmail.com}

\begin{abstract}
We prove that every finite simple group of Lie type $G$ can be generated by three regular unipotent elements. 
In certain cases we show that two regular unipotents are sufficient to generate $G$.
\end{abstract}

\subjclass[2000]{20F05, 20D06} 
\keywords{Simple groups of Lie type; regular unipotent element;  generation}
\maketitle


\section{Introduction}

Problems of the generation of finite simple groups by certain their elements 
attract a  significant attention of group theorists and have many applications. 
For many occasions one requires a group to be generated by specific elements. This kind of generation problems forms a large area of research in group theory. One of the aspects of this research is obtaining a good upper bound for the minimal number of conjugates of a given group element to generate the group in question.    

The most universal results have been obtained by  Guralnick and Saxl \cite{GS}.
For a finite nonabelian simple group $G$ and 
$g\in G$, let ${\rm gn}_G(g)$ be the minimal number of conjugates of $g$ that generate $G$, and let 
${\rm gn}(G)=\max_{1\neq g\in G}{\rm gn}_G(g)$.
In \cite{GS} the authors obtained a rather sharp lower bound for ${\rm gn}(G)$.
However, for a specific $g\in G$ the actual minimal number of conjugates of $g$ needed to generate $G$ dramatically depends on the choice of $g$.
The problem of understanding the dependence of ${\rm gn}_G(g)$  on the choice of $g\in G$ is of great interest. Some results in this direction are proved in the literature (for instance, see \cite{RZ}  that 
provides optimal bounds for almost all sporadic groups), but in general the problem is far from a satisfactory solution. A natural
and important question is that of determining $g\in G$ with ${\rm gn}_G(g)=2$. 
Examples of such elements can be found in various publications.  

L\"ubeck and Malle proved in \cite{LM} that every finite simple group of exceptional Lie type (apart from the Suzuki groups) is generated by two elements of respective order $2$ and $3$. 
Similar results for classical groups are obtained in \cite{P,PT,PT2,PT3}. It is shown in \cite{Sn} that, for $1\neq u\in G$ unipotent, $G=\langle u,h\rangle$ for a suitable  semisimple element $h\in G$.  In \cite{VN} the authors prove (modulo some open conjectures) 
that all but two nonabelian finite simple groups are generated by three conjugate involutions.

Our interest in this paper concerns regular unipotent elements of finite simple groups of Lie type.  
We conjecture that if $G\neq PSL_2(q)$ with $q$ even, then $G$  is generated by two regular unipotent elements that can be  chosen conjugate in $G$ provided that $G\neq PSL_2(9)$, see below for details.
\smallskip

There is a certain connection of these problems with those on presentations of an arbitrary group element as a product of two elements of special kind. An important contribution was done by Gow \cite{Gow} who proved that every semisimple non-identity element of a finite simple group $G$ of Lie type is a product  of two  elements from an arbitrary conjugacy class of regular semisimple elements. This leads to the conclusion that $G$ is generated by three  elements from any conjugacy class of regular semisimple elements \cite[Lemma 2.8]{DZ}.
Ellers, Gordeev and Herzog \cite[Theorem H, p. 344]{EGH}  proved that every non-identity element of $G$ is a product of two unipotent elements (take there $h=1$). The argument \cite[Lemma 2.8$(2)$]{DZ} implies that  $G$ is generated by three unipotent elements, one can be chosen arbitrarily.

\begin{conje}\label{conj1}
Let $G$ be a quasisimple finite group of Lie type such that    
$G\neq SL_2(q)$ with $q$ even. Then $G$ is generated by two regular unipotent elements. \end{conje}
 
Groups $SL_2(q)$  with $q>2$ even are not generated by two unipotent elements as these are involutions, and two involutions generate either an abelian group or a dihedral one. 

Recall that the regular unipotent elements of $G$ are not always conjugate. 
A  stronger version of Conjecture \ref{conj1} is that  $G$ is generated by two conjugate regular unipotents   if, additionally, $G\neq PSL_2(9)$. 
If $G= PSL_2(q)$ with $q\geq 5$ odd, then the validity of the stronger version of Conjecture \ref{conj1} is contained in \cite[Lemma 3.1]{GS}.
If $p$ is a good prime for $\GG$ then the regular unipotents are conjugate in $\Aut(G)$,
see for instance \cite[Lemma 5.2]{TZ04}.  (Here $\GG$ is a simple algebraic group used to define $G$, see below. Recall that a prime $p$ is good unless $p=2$ and $\GG$ is not of type $A_n$, $p=3$ for $\GG$ of the exceptional Lie type, and $p\leq 5$ for $\GG$ of type $E_8$. In addition, if $G$ is of adjoint type then $Z(G)=1$.)
This implies (for $p$ good) that, if $G$ is generated by $k$ conjugate regular unipotents, then for every regular unipotent element $u\in G$ there are $k$ conjugates of $u$ that generate $G$.

Our main results are the following.
\begin{theo}\label{qs3}
Every finite quasisimple group of Lie type is generated by three regular unipotent elements.\end{theo}

\begin{theo}\label{t22}
Conjecture {\rm \ref{conj1}} holds for the following groups of Lie type:
\begin{enumerate}[$(1)$]
\item $SL_2(q)$ with $q\geq 4$, and  $SL_3(q)$;
\item $SU_3(q)$ with $q\geq 3$;
\item $Sp_4(q)$ and $G_2(q)$ with $q$ odd;
\item ${}^2B_2(2^{2m+1})$, ${}^2G_2(3^{2m+1})$ and  ${}^2F_4(2^{2m+1})$, with  $m\geq 1$;
\item ${}^3D_4(q)$;
\item $SU_4(2)$,  $SU_4(3)$, $SU_5(2)$, $SU_5(3)$, $Sp_6(9)$, $G_2(4)$, $F_4(2)$, ${}^2F_4(2)'$ and ${}^2E_6(2)$;
\item  $SL_n(q)$ with  $n\geq 4$. 
\end{enumerate}
Furthermore, in cases $(1)$--$(6)$, and in  $(7)$ for $q$  a prime, the stronger version of Conjecture {\rm \ref{conj1}} holds.
\end{theo}

For other quasisimple  groups of Lie type we prove  a slightly weaker bound.

\begin{theo}\label{t33}
The following finite quasisimple groups of Lie type are generated by three conjugate regular unipotent elements:
\begin{enumerate}[$(1)$]
\item $SL_2(q)$ with $q\geq 4$ even, $SL_2(9)$ and $SL_n(q)$ with $n\geq 3$;
\item $SU_n(q)$ with $n\geq 3$ odd;
\item $Sp_{2n}(q)$ with  $n\geq 2$;
\item $\Omega_{2n+1}(q)$ with $n\geq 3$ and $q$ odd;
\item $\Omega_{2n}^+(q)$ with $n\geq 4$ and $q$ even;
\item $\Omega_{2n}^-(q)$ with $n\geq 4$;
\item $SU_4(q)$,  with $q$ even;
\item $G_2(q)$, $F_4(q)$ and  ${}^2E_6(q)$.
\end{enumerate}
\end{theo}

For a subgroup $X$ of  a finite group $G$ of Lie type, let $X_{\un}$ be the subgroup generated by the unipotent elements of $X$.  
Our strategy in proving Theorem \ref{t22} is in showing that there is a maximal subgroup $X$
of $G$ such that $X_{\un}$ is  generated by two regular unipotents. 
In many occasions for this a maximal parabolic subgroup can be taken. So we state our result in this direction as follows.

\begin{theo}\label{p33}  
Let $G$ be one of the following quasisimple groups:
$SL_n(q)$ where $n\geq 3$, 
$SU_n(q)$ where $n\geq 6$ is even,
$Sp_{2n}(q)$ where $n\geq 3$, 
$\Omega_{2n}^\pm(q)$ where $n\geq 4$, 
$\Omega_{2n+1}(q)$ where $n\geq 3$ and $q$ is odd,
${}^2E_6(q)$ where $q$ is odd,
$E_6(q)$, $E_7(q)$, $E_8(q)$.
Then $G$ contains a maximal parabolic subgroup
$P$ such that $P_\un$ is generated by two regular unipotent elements. 
\end{theo}

\subsubsection*{Notation}
We denote by $\FF_q$ the finite field of $q$ elements. In the following, $J_n$ denotes an upper triangular unipotent Jordan block of size $n$. We write $\diag(d_1\ld d_k)$ for a block-diagonal matrix with diagonal blocks $d_1\ld d_k$. We denote by ${}^tA$ the transpose of a matrix $A$.
For a group $G$ we write $Z(G)$ for the center of $G$, $G'$ for the derived subgroup and $O_p(G)$ for the maximal normal $p$-subgroup of $G$ (where $p$ is a prime).
For subgroups $X,Y\leq G$ we write $[X,Y]$ for the group generated by $xyx\up y\up$ with $x\in X,y\in Y$. For $g\in G$, $|g|$ denotes the order of $g$.

Let $\GG$ be a connected reductive algebraic group in characteristic $p$. We view it as $\GG(F)$, where $F$ is an arbitrary \acf of characteristic $p$, since our reasoning do not depend on the choice of $F$.   Let $\si$ be a Frobenius endomorphism of $\mathbf{G}$,   and let  $G=\GG^\si$ be a finite group of Lie type. 
Then $\si$ stabilizes some Borel subgroup $\mathbf{B}$ of $\mathbf{G}$ and a maximal torus $\mathbf{T}$ of $\mathbf{B}$.  
These determine  a root system of $\mathbf{G}$ and, for every root $\al$, the root subgroups $x_{\al}$, each  is isomorphic to the additive group of the ground field $F$. By $\Phi(G)$ we denote  the set of positive roots and $S$ the subset of simple roots. Every unipotent element of $\mathbf{B}$ can be written as the 
product $\Pi _{\al\in \Phi}x_\al (t_\al)$ for some $t_\al\in F$. In these terms, a regular unipotent element of $\mathbf{B}$ can be defined as such that 
$t_\al\neq 0$ for every $\al\in S$, and a regular unipotent element of 
$\mathbf{G}$ can be defined as a conjugate (in $\mathbf{G}$) of a regular unipotent element of 
$\mathbf{B}$ (see   \cite[Proposition~5.1.3]{CarF}). We keep this term for elements of $G=\GG^\si$ that are regular unipotents in $\GG$, and also apply this to the group ${}^2 F_4(2)'$ which is not of such a form.

Recall that parabolic subgroups $\mathbf{P}$ of $ \mathbf{G}$ are those containing 
a Borel subgroup of $\GG$, and those containing our fixed subgroup $\mathbf{B}$ are referred here as standard parabolics.   
Every parabolic subgroup $\mathbf{P}$ is conjugate to a standard  one,  
which, in turn,  is a semidirect product $\mathbf{P} = \mathbf{U}\mathbf{L}$, where  $\mathbf{L}$ is a  connected reductive group containing $\mathbf{T}$ and $\mathbf{U} $ is the unipotent radical of $\mathbf{P}$. 

By general theory, we can choose $\mathbf{B}$ and $\mathbf{T}$ $\si$-stable.
So $\mathbf{T}^\si <  \mathbf{B}^\si < G$. If $\mathbf{P}  $ is a standard parabolic and $\si(\mathbf{P} )=\mathbf{P} $, the group $P:=\mathbf{P}^\si$ is called a standard parabolic subgroup of $G$. (For our purposes we can always assume that $P$ is standard.)  In addition, if $\si(\mathbf{P})=\mathbf{P}$ then there exists a  Levi subgroup $\mathbf{L}$ of $ \mathbf{P}$ such that $\si(\mathbf{L})=\mathbf{L}$.

Let $N=N_{\mathbf{G}}(G)$. The automorphisms of $G$ induced by inner automorphisms of $N$ are called diagonal.  

For an algebraic group $\mathbf{G}$, or a finite group $G$ of Lie type, we denote by $ \mathbf{G}_{\un}$ and $G_{\un}$ the subgroups generated by the unipotent elements. 
Every connected reductive algebraic group $\mathbf{G}$ is a central product $Z(\mathbf{G})\cdot \GG_{\un}$, where $\GG_{\un}$ is a semisimple algebraic group. 
In particular, $\mathbf{L}=Z(\mathbf{L})\cdot\mathbf{L}_{\un}=\mathbf{T}\cdot \mathbf{L}_{\un}$, and hence the automorphisms of $\mathbf{L}_{\un}^\si$ arising from the action of $ \mathbf{T}^\si$ on 
$\mathbf{L}_{\un}^\si$ are diagonal.

\section{Preliminaries}

\bl{qg1}\cite[Lemma 2.6]{tz13}
Let $\GG$ be a reductive algebraic group,  and $\mathbf{P}$  a parabolic subgroup of $\GG$ properly containing a  Borel subgroup of $\GG$. 
Let $\mathbf{L}< \mathbf{P}$ be a Levi subgroup of $\mathbf{P}$.
Let $u\in \mathbf{P}$ be a regular unipotent. Then $u\notin \mathbf{L}$; in addition, the projection of $u$ into $\mathbf{L}$ is a regular unipotent in $\mathbf{L}$. 
Consequently, every regular unipotent  $u\in \mathbf{L}$ is the projection of some regular unipotent element of~$\mathbf{P}$.
\el

The last claim is not stated in \cite[Lemma 2.6]{tz13}, but it follows from the previous claim as the regular unipotent elements of an algebraic group are 
conjugate.

Now, let $\si$ be a Frobenius endomorphism of $\mathbf{G}$, a simple algebraic group, and let  $G=\GG^\si$ be a finite group of Lie type. 
Then there are $\si$-stable conjugates of $\mathbf{P}$ and $\mathbf{L}< \mathbf{P}$. 
Let $P=\mathbf{P}^\si$ and $L=\mathbf{L}^\si$. Then  $P$ is a parabolic subgroup of $G$ (in sense of the BN-pairs theory), and $L$ is a Levi subgroup of $P$, that is, a complement of $O_p(P)$ in $P$. In this notation we have:

\bl{n26} 
The projection of a regular unipotent element $u\in P$ is a regular unipotent of $L$. 
Moreover, every regular unipotent element of $L$ is a projection of a regular unipotent of $P$ into $L$.\el

\bp The first assertion follows from Lemma \ref{qg1} as the notion of a regular unipotent element of $G$ and $L$ is described in terms of algebraic groups $\GG$ and $\mathbf{L}$ \cite[Proposition 5.1.3]{CarF}. Recall that every unipotent element  $u$ in a standard Borel subgroup is expressed as a product 
$\Pi_{\al}(x_{\al}(t_\al))$ with the product ranging over positive roots of $G$, and $u$   is regular \ii $t_\al\neq 0$ for all simple roots $\al$. We say that  $\mathbf{L}$ is standard if the simple roots of $\mathbf{L}$ are simple roots of $\mathbf{G}$. 
If $\mathbf{L}$ and $\mathbf{P}$ are $\si$-stable and standard, we call $L=\mathbf{L}^\si$ and $P=\mathbf{P}^\si$ the standard parabolic and Levi subgroups respectively. 

If $G$ is not of twisted type then the second claim follows by exactly  the same argument as in the proof of \cite[Lemma 2.6]{tz13}. 
Suppose that $G$ is twisted.
Then the argument requires some adjustment. Note that standard $\si$-stable parabolic subgroups of $G$ are determined by some $\si$-orbits on the set of 
simple roots of $\mathbf{G}$ and on the respective sets of root subgroups $x_\al(t_\al)$, see \cite[Proposition 26.1]{MTe}. 
Unipotent elements in a standard Levi subgroup are  products $\Pi_{\al}(x_{\al}(t_\al))$ over the positive roots of 
$\mathbf{L}$, which are  $\si$-stable in the sense that each  term $x_{\al}(t_\al)$ occurs together with 
$x_{\si(\al)}(\si(t_\al))$. Note that $t_\al\neq 0$ implies $\si(t_\al)\neq 0$. 
So, $u\in L$ is regular unipotent if all $t_\al\neq 0$. If $S'$ is the set of simple roots of $G$ that are not simple roots of  $\mathbf{L}$, then $\si(S')=S'$.  Then we can extend the above product by the multiple $\Pi _{\al\in S'}x_{\al}(1)$ to obtain a regular unipotent element $u'\in P$ that is regular in   $G$ and whose projection into  $\mathbf{L}$ coincides with $u$. 
\enp

Note that Levi subgroups of $P$ are not always conjugate but it suffices for our purpose to deal with the standard Levi subgroup.

\bl{bor} 
Let $G$ be a simple group of Lie type in defining characteristic $p>0$ and $U$ be a \syl of $G$. Suppose that $H$ is a maximal subgroup of $G$ containing $U$.
Then $H$ is a maximal parabolic subgroup of $G$. In addition, 
$H$ is the only maximal subgroup of $G$ containing $H_{\un}$. 
\el

\bp It is well known that $G$ is group with BN-pair \cite[Theorem 24.10]{MTe}. 
These groups satisfy the assumption of  \cite[Ch. IV, \S 2.7, Theorem 5]{Bo}, which states that $H$ is either parabolic or normal in $G$, whence the result. 
For the additional claim, if $H$ is a maximal parabolic, then $H=N_G(O_p(H))$.
Suppose that $Q$ is a  maximal subgroup of $G$ containing $H_\un$. Then $H\leq Q$ by a theorem of Borel and Tits, see \cite[13.1]{GL}.
As $H$ is maximal, we have $H=Q$.
\enp

\bl{nn6} 
Let $G$ be a quasisimple group of Lie type,  $P$  a maximal parabolic of $G$ and let $u\notin P$ be a unipotent element. 
Then $\lan u,P_{\un}\ran=G$. \el

\bp 
Set $X=\lan u, P_\un \ran$. By Lemma \ref{bor}, $P$ is the only maximal subgroup of $G$ containing $P_\un$. Since $u \not \in P$,  $X$ cannot be contained in $P$ and hence $X=G$.
\enp
  
The \f lemma can be useful.

\bl{md8} 
Let $G < GL_n(q)=GL(V)$, and let $g,h\in G$. Suppose that the Jordan form of $g$ and $h$ is $J_n$. 
Suppose that, with respect to some basis of $V$, the matrix of $g$ is upper triangular and the matrix of $h$ is lower triangular. Then $X=\lan g,h\ran$ is \irt \el

\bp There is a unique series $\{0\}=V_0\subset V_1 \subset\ldots \subset V_{n}=V$  of $g$-invariant subspaces of $V$ with  $\dim V_i=i$ for $i=1\ld n$. Let $v_i\in V_i\setminus V_{i-1}$
for $i=1 \ld n$, and set $V_i'=\lan  v_{i+1}\ld v_{n}\ran $. As $h$ is lower triangular and unipotent, $hV_i'=V_i'$ and $(h-\Id)v_i\in V_i'$.  

Suppose that $X$ is reducible. Then  $hV_i=V_i$ for some $i\in\{1\ld n-1\}$ and hence $hv_i\in V_i$. This is a contradiction, unless $hv_i=v_i$. 
Since $h v_n=v_n$, we get the absurd $\dim C_V(h)>1$.\enp

\bl{ss42} The groups $SL_4(2)$ and  $Sp_4(2)$ are generated by two conjugate elements with Jordan form~$J_4$.\el

\bp Observe that $SL_4(2)\cong Alt(8)$ and  $Sp_4(2)\cong Sym(6)$. Let $G=SL_4(2)$. 
Then $C_G(J_4)$ is an abelian $2$-group of order $8$. 
So, $J_4$ lies is in class $4B$ in the notation of \cite{at} and is realized as the permutation $(1,2,3,4)(5,6)$ in $Alt(8)$. 
Let $t= (2,5,8,4,6,3,7)$. Then $h:= t\up g t =(1,7,6,8)(2,4)$ and  $g^2h  =(1,7,6,8,3) $ is of order $5$. The group $\lan g,h\ran $ is primitive and contained in $Alt(8)$. 
By \cite[Theorem 3.3E]{DiM},  $\lan g,h  \ran =Alt(8)$.

Let $G=Sp_4(2)\cong Sym(6)$. The isomorphism $Sym(6)\stackrel{\cong}{\longrightarrow} Sp_4(2)$ can be chosen from taking the non-trivial \ir constituent of the natural permutational module of degree $6$ over $\FF_2$. 
Then we choose for $g$ the $4$-cycle  $(1,2,3,4)$ in $Sym(6)$,
by viewing $Sym(6)$ as the group of permutations on $\{1,2,3,4,5,6\}$.
Let $t= (2,6,4,5,3)$. Then $h:= t\up g t =(1,3,5,6)$ and 
$(g^2h)^2 = (3,6,5)$ is of order $3$. The group $\lan g,h\ran $ is primitive and not 
contained in $Alt(6)$. By \cite[Theorem 3.3E]{DiM},  $\lan g,h\ran =Sym(6)$. \enp
 
We describe the cases where the regular unipotent elements of $G$ are conjugate.
By \cite[Lemma 5.2]{TZ04}, this is the case \ii  $|Z(G)|=1$ and $p$ is a good prime for $\GG$  (as defined prior to Theorem \ref{qs3} above).
If $p$ is good and $Z(G)=1$ then  either $G=SL_n(q)$ with $\gcd(q,n-1)=1$ or $G=F_4(q),G_2(q)$ with $p>3$ or $E_8(q)$ with $p>5$.

\section{Groups $PSL_n(q)$}

\bl{gt1} If $G=SL_n(q)$ with $n>2$,  then  Conjecture {\rm \ref{conj1}} is true. \el  

\bp If $q$ is a prime then the result is contained in  \cite{GT1}, unless $(n,q)=(4,2)$.  For  $G=SL_4(2)$ see Lemma \ref{ss42}. So we can assume that $q>3$.
 
By \cite[Theorem 6]{Lev1}, if $q>3$ then every non-central element of $G$  is a product of two $GL_n(q)$-conjugates of $J_n$; these are of course regular unipotents in $G$. In particular, an \ir element $t$ of order $(q^n-1)/(q-1)$ can be written as a product $t=gh$, where $g,h$ are regular unipotents in $G$.
So $\lan g,h\ran$ is an \ir subgroup of $G$ containing  
$T=\lan t\ran$. It is well known that $T$ acts transitively on the lines of $V$, the underlying space of $G$, and, consequently, irreducibly on $V$. Groups $H$ acting transitively on the lines of $V$ are determined in  \cite{he1}, see also \cite[p. 512--513]{Lie}. We use the list provided in \cite{Lie} to exclude those containing no element of order 
$(q^n-1)/(q-1)$ and no regular unipotent. 
Let $H'$ be the derived subgroup of $H$ (the case $H=H'$ is not excluded).

Suppose that $H\neq SL_n(q)$. By \cite[p. 512--513]{Lie}, one of the \f holds:
 
\begin{itemize}
\item[(A)] $H'\cong SL_{n/k}(q^k)$;
\item[(B)] $H'= Sp_n(q)$, $n\geq 2$ even, or $G_2(q)'$ for $n=6$ and $q$ even;
\item[(C)] $H$ normalizes an \ir extraspecial $2$-group $E$, and either $E\cong Q_8$, $n=2$ and $q\in\{5,7,11,23\}$ or  $n=4,q=3$ and $E\cong Q_8\circ D_8$ and $H'/E \leq S_5$;
\item[(D)] either $H'$ normalizes $SL_2(5)$, $n=2$ and $q\in\{9,11,19,29 ,59\}$ or $H'=A_6,A_7 < SL_4(2)$, or $H'=SL_2(13)< Sp_6(3)$.
\end{itemize}

Suppose that $g,h\in H$. Then $H$ is not a group from (A), as such a group  contains no  element with Jordan form $J_n$.

If (B) holds then $n>2$ is even,   and the order of an \ir cyclic subgroup of $H'$ divides $q^{n/2}-1 $ or $q^{n/2}+1 $. The maximum order of an element of $PCSp_{n}(q)$ is less than $q^{\frac{n+2}{2}}/(q-1)$ (see \cite[Lemma 2.10]{GMPS}).
Since $q^{\frac{n+2}{2}} < q^{n}-1$, we obtain that $t\not \in H$.
If $n=6$ and   $H'=G_2(q)'$ then $H' <  Sp_6(q)$, so this case is ruled out.

As $n>2$ and $q$ is not a prime,  cases (C) and (D) cannot occur. This completes the proof.  \enp

\begin{corol}\label{c21} 
Let $G=SL_n(q)$, where $n>2$ and $\gcd(n,q-1)=1$.  
Then  $G$ is generated by two conjugate regular unipotent elements. 
\end{corol}

\bp 
The result follows from Lemma \ref{gt1}, as in this case the regular unipotent elements are conjugate in $G$. Indeed, let $d$ be the number of conjugacy classes of regular unipotents in $G$. Recall that $d$ equals $|Z( G)|$ (for $G=SL_n(q)$), see \cite[Lemma 5.2(i)]{TZ04}.  
\enp

\begin{rem}\label{SL29}
The case $n=2$ is treated in \cite[Lemma 3.1]{GS}, assuming $q\geq 4$. In particular, $SL_2(q)$ for $q\neq 9$ odd is generated by two conjugate regular unipotents,
while $SL_2(9)$ is generated by three conjugate regular unipotents.
One can also verify that $SL_2(9)$ is generated by the regular unipotents $a=J_2$ and $b=\left(\begin{smallmatrix} 1 & 0 \\
\beta & 1\end{smallmatrix}\right)$, where $\beta \in \FF_9$ is such that $\beta^2=\beta+1$.
In fact, $ab^2$ and $[a,b]$ have respective order $8$ and $10$, but no maximal subgroup of $SL_2(9)$ contains elements of both such orders.
\end{rem}

\bl{t2r} 
The matrix  ${}^t J_n$ is conjugate in $SL_n(q)$ to $J_n$ or to $J_n\up$.
\el

\bp Consider the matrix $M$ with $1$ at the second diagonal and $0$ elsewhere. Then $M J_n M\up={}^t J_n$. If $q=2$ then $\det M=1$ and the result follows.  If $q$ is  odd  then $\det M=1$ \ii $n\equiv 0,1 \pmod 4$. Otherwise, $\det M=-1$. If $n\equiv 3 \pmod 4$ then $-M\in SL_n(q)$, and the result follows. 

Let $n\equiv 2 \pmod 4$. If $-1=\lam^2$ for some $\lam\in \FF_q$ then $\det(\lam\Id)=-1$ and $\lam M\in SL_n(q)$. 
Note that $J_n$ is not conjugate to ${}^tJ_n$. Indeed, if $xJ_nx\up={}^tJ_n$ then $Mx\in C_{GL_n(q)}(J_n)$. As $ C_{GL_n(q)}(J_n)$ consists of the upper triangular matrices with scalar diagonal, $\det(Mx)=1$ implies the existence of a scalar matrix with determinant $-1$, which is false.

We show that ${}^tJ_n$ is conjugate in $SL_n(q)$ to $J_n\up$. Set $N=\Id -(J_n-\Id)$. Then $N$ is a regular unipotent with   $-1$ above the diagonal. 
Let $M_1$ be the antidiagonal matrix $\mathrm{antidiag}(+1,-1,\ldots,+1,$ $-1)$.
Then $M_1{}^tJ_nM_1\up=N$ and $\det M_1=1$. Note that $J_n\up$ is another matrix with  $-1$ above the diagonal. By   \cite[Lemma 5.1]{TZ04}, $J_n\up$ and $N$ are conjugate in  $  B$, the group of  upper triangular  matrices in $G$.
\enp

\begin{corol}\label{co1} 
Let $G=SL_n(q)$, where $n\geq 2$ and $q$ is a prime. 
Then  $G$ is generated by two conjugate regular unipotent elements. 
\end{corol}

\bp  
If $(n,q)=(4,2)$ then the result follows from Lemma \ref{ss42}.  
Otherwise, by \cite{GT1}, $G=\lan J_n, {}^tJ_n\ran$. 
By Lemma \ref{t2r}, the matrix ${}^tJ_n$ is conjugate in $SL_n(q)$ to $J_n$ or to $J_n^{-1}$. Since $\lan J_n, {}^tJ_n\ran=\lan J\up_n,{}^t J_n\ran$, 
the result follows. 
\enp

From the previous corollary, it follows in particular that the group $SL_2(3)$ is generated by two conjugate regular unipotents. This fact will be used in Lemmas \ref{ee6} and \ref{2e6}.

\bl{cc3} 
Suppose that the cyclic groups generated by the regular unipotents are conjugate. If $G$ is generated by $k$ regular unipotents, then $G$ is 
generated by $k$ conjugate regular unipotents.  
\el

\bp 
Note that $u$ and $u^m\neq 1$ generate the same cyclic subgroup of $G$. By assumption, $G$  is generated by $k$ regular unipotents $u_1\ld u_k$, say. \itf $G$  is generated by $k$ subgroups  $\lan u_1\ran\ld \lan u_k\ran$. As every subgroup $\lan u_i\ran$ is conjugate to $\lan u\ran$, the result follows.
\enp

\begin{propo}\label{zg1}
Let $G=SL_n(q)$, where $n\equiv 2 \pmod 4$ and $\gcd(q-1,n)=2$. Then  $G$ is generated by two conjugate regular unipotent elements. 
\end{propo}
	
\bp 
By \cite[Lemma 5.2(i)]{TZ04}, the number of conjugacy classes of regular unipotents in $G$ is $|Z(G)|=2$.  Let $u\in G$ be a regular unipotent. 
By \cite[Theorem 1.7(ii)]{TZ04}, $u$ is not rational if $q$ is odd, $n=2m$ is even and $n/\gcd(n,q-1)=2m/|Z(G)|=m$ is odd. \itf each conjugacy class of regular unipotents in $G$ meets  the cyclic group $\lan u\ran$.  Therefore,
if $G=\lan u,v\ran$ for some regular unipotent $v\in G$, by replacing $v$ by a suitable power of $v$, we obtain two conjugates of $u$ that generate $G$.
\enp

Let $d$ be the number of conjugacy classes of regular unipotents in a finite simple group $G$ of Lie type. Recall that $d$ equals  $|Z(G)|$ 
if $G=\GG^\si$ for some Steinberg endomorphism $\si$ of $\GG$, see \cite[Lemma 4.7(ii)]{TZ04}.
If $G= Spin_{2n+1}(q)$ with $q$ odd then $|Z(G)|=2$, ($p$ is good) and 
regular unipotents are not conjugate \ii  $2n+1\neq  \pm 1\pmod 8$, see \cite[Theorem 1.7(iv)]{TZ04}. 
If $u$ is regular unipotent then the generators of $\lan u \ran$ form (at most) two conjugacy classes by \cite[Theorem 1.5]{TZ04}. 
So the above reasoning works and in these cases $G$ is generated by three conjugate regular unipotents.

If every regular unipotent is rational then $k$ coincides with the number of cyclic subgroups generated by  regular unipotents. These cases are listed in
\cite[Theorem~1.8]{TZ04} for $G\neq E_8(q)$. In particular, if $q$ is even or a square then this is the case, so Lemma \ref{cc3} is not useful if we wish to   improve the result of Lemma \ref{gt1} by stating that $G$ is generated by two conjugate unipotents.

\section{Generation by two conjugate regular unipotent elements}

\begin{corol}\label{d66} 
The group $\Omega_6^+(q)$ is generated by two regular unipotent elements.
\end{corol}

\bp Note that $\Omega_6^+(q)$ is isomorphic to $PSL_4(q)$ and the isomorphism
in question arises from an \irr $\phi:SL_4(q)\ra GL_6(q)$. Then the regular unipotent $u\in SL_4(q)$ remains a  regular unipotent in  $\Omega_6^+(q)$. The Jordan form of it is 
$\diag(J_2,J_4)$ if $q$ is even, and it is $\diag(J_5,J_1)$ if $q$ is odd.
\enp

 \bl{sz2g}
For every regular unipotent element $u$ of $G=Sz(2^{2m+1})$, $m\geq 1$, there exists $g \in G$ such that $G=\lan u,u^g\ran$.
\el

\bp 
Let $u$ be a regular unipotent element of $G=Sz(2^{2m+1})$. Since $G\leq Sp_4(q)$,  the Jordan form of $u$ is $J_4$, whence $|u|=4$.
As shown in \cite{Suz}, $G$ admits a unique class of involutions and exactly two classes of elements of order $4$, respectively those of $u^2$, $u$ and $u^{-1}$.
Always by \cite{Suz}, $G$ can be generated by an element $x$ of order $2$ and an element $y$ of order $4$.
Then $y$ is regular unipotent and $H=\langle y,y^x\rangle$ is a normal subgroup of $\langle x, y \rangle=G$. Since $G$ is simple, it follows that $H=G$. 
\enp

The knowledge of the character table of a group $G$ (and of the list of its maximal subgroups) can be very useful to our purposes.
Namely, given a finite group $G$, let $c_1,c_2,c_3$ be (not-necessarily distinct) conjugacy classes of $G$.
Denote by $\Delta_G(c_1,c_2,c_3)$ the number of distinct triples $(g_1,g_2,g_3)$, where $g_3$ is a fixed element of the class $c_3$
and $g_1 \in c_1$, $g_2\in c_2$ are such that $g_1g_2=g_3$.
This structure constant can be computed using the (complex) character table of $G$. Namely, it is given by the formula
\begin{equation}\label{Delta}
\Delta_G(c_1,c_2,c_3)=\frac{|c_1|\cdot |c_2|}{|G|}\cdot \sum_{i=1}^r\frac{\chi_i(g_1)\chi_i(g_2)\overline{\chi_i(g_3)} }{\chi_i(1)},
\end{equation}
where $\chi_1,\ldots,\chi_r$ are the irreducible complex characters of $G$.

\bl{U5(2)}
For every regular unipotent $u$ of $G=SU_5(2)$, there exists $g \in G$ such that $G=\lan u,u^g\ran$.
\el

\bp
The regular unipotent elements of $G=SU_5(2)$ have Jordan form $J_5$ and belong to the class $8a$. 
Using the character table of $G$, one can compute  $\Delta_G(8a,8a,11a)=53416$.
So, given a regular unipotent element $u$, there exists $g \in G$ such that $uu^g$ has order $11$. Now, the only maximal subgroup of $G$ whose order is divisible by $11$ is $PSL_2(11)$, but this subgroup
does not contain elements of order $8$. Hence, $\lan u,u^g\ran=G$.
\enp

\bl{2E6(2)} 
For every regular unipotent $u$ of $G={}^2E_6(2)$, there exists $g \in G$ such that $G=\lan u,u^g\ran$.
\el

\bp
The group $G={}^2E_6(2)$ has three classes of regular unipotent elements. In fact, $G$ has nine classes of elements of order $16$, three of them ($16a$, $16b$ and $16c$) with centralizer of size $256$, the other six with centralizer of size $128$. 

Let $c$ be any conjugacy class of regular unipotents and let $c_3$ be the class $19a$.
Since $\Delta_G(c,c,c_3)>0$, for every regular unipotent element $u\in G$ there exists $g\in G$ such that the product $u u^g$ has order $19$.
The list of the maximal subgroups of $G$ is provided  in \cite{at,C23}.
It is easy to see that the only maximal subgroup of $G$ whose order is divisible by $19$ is $PSU_3(8):3_1$. However, this subgroup does not contain elements of order $16$, proving that $G=\lan u, u^g\ran$.\enp

\bl{2g2}
For every regular unipotent $u$ of $G={}^2G_2(3^{2m+1})$, $m\geq 1$,  there exists $g \in G$ such that $G=\lan u,u^g\ran$.
\el 

\bp
The maximal subgroups of $G={}^2G_2(3^{2m+1})$ are described in \cite{BHRD, Kl88b}. In particular, as shown in \cite{We}, the only proper overgroup of
the cyclic torus $T$ of order $t=3^{2m+1}+3^{m+1}+1$ is its normalizer $H=T:6$.
Note that the unipotent regular elements of $G$ have order $9$, as $G\leq \Omega_7(3^{2m+1})$.

Our aim is to show that there exist two conjugate regular unipotents whose product has order $t$.
This goal can be achieved by computing the structure constants for $G$ thanks to \eqref{Delta}. In particular, we make use of \cite{Chevie}, that contains the character table of $G$ and all the relevant information we need. So, following the notation of \cite{Chevie}, let $c$ be a conjugacy class consisting of regular unipotents (so,
$j=5,6,7$) and let $c_{14}$ be a class of elements of order $t$ (with the parameter  $\mathtt{iI}=1$).
Computing $\Delta_G(c_j,c_j,c_{14})$ we obtain the value
$$3^{10m+3} - 3^{8m+2} + 4\cdot 3^{6m+1} + 2\cdot 3^{7m+2}  + 2\cdot   3^{5m+1 } + 2\cdot  3^{4m+1} +  3^{3m+1}$$
for $j=5$, and
$$3^{10m+3}  - 3^{8m+2} - 3^{7m+2}  -2\cdot  3^{6m+1} - 3^{5m+1}$$
for $j=6,7$. Note that $\Delta_G(c_j,c_j,c_{14})>0$ for all $j=5,6,7$. This implies that, given a regular unipotent element $u$, there exists $g \in G$ such that $u u^g$ has order $t$. We conclude that $\lan u,u^g\ran $ contains $T$, whence the result.
\enp

\bl{3d4q}
For every regular unipotent $u$ of $G={}^3D_4(q)$,  there exists $g \in G$ such that $G=\lan u,u^g\ran$.
\el 

\bp
The maximal subgroups of $G={}^3D_4(q)$ are described in \cite{Kl88a, BHRD}. In particular, as shown in \cite{We}, the only proper overgroup of
the cyclic torus $T$ of order $t=q^4-q^2+1$ is its normalizer $H=T: 4$.
Note that the regular  unipotent elements of $G\leq \Omega_8^+(q^3)$ have Jordan form $\diag(J_2,J_6)$, so they have order $8$ when $q$ is even.

We proceed, as done in proving Lemma \ref{2g2}, with using the character table of $G$ in order to compute the structure constants of $G$.
Following the notation of \cite{Chevie}, let $c_j$ be a conjugacy class consisting of regular unipotents 
and let $c_{k}$ be a class of elements of order $t$ (with the parameter  $\mathtt{kK}=1$).
So, $j=7,8$ and $k=28$ if $q$ is even,  $j=7$ and $k= 32$ if $q$ is odd. Computing $\Delta_G(c_j,c_j,c_{k})$ we obtain the value 
$$\frac{1}{4}q^{20}-\frac{1}{4} q^{18} + \frac{1}{4}q^{16} -\frac{1}{2}q^{14} +q^{12} -q^{10} +\frac{1}{2} q^8$$
if $q$ is even, and
$$q^{20}-q^{18}+q^{16}-2q^{14}+2q^{12}-2q^{10}$$
if $q$ is odd.
Note that $\Delta_G(c_j,c_j,c_{k})>0$ for all values of $j$ and $k$. This implies that, given a regular unipotent element $u$, there exists $g \in G$ such that $u u^g$ has order $t$. We conclude that $\lan u,u^g\ran $ contains $T$, whence the result.
\enp

\begin{lemma}\label{Ree}
For every regular unipotent element $u$ of $G={}^2F_4(2^{2m+1})$, $m\geq 1$, there exists $g \in G$ such that $G=\lan u,u^g\ran$.
\end{lemma}

\begin{proof}
Set $q^2=2^{2m+1}$. The maximal subgroups of $G={}^2F_4(2^{2m+1})$ are described in \cite{M}. 
In particular, the only proper overgroup of
the cyclic torus $T$ of order $t=q^4+\sqrt{2} q^3+q^2+\sqrt{2} q+1 $ is its normalizer $H=T: 12$.
Note that the regular unipotent elements of $G$ have order $16$. So, to prove the result it suffices to show that, given a regular unipotent element 
$u \in G$, there exists $g \in G$ such that the product $uu^g$ has order $t$.

We proceed as done in proving Lemma \ref{2g2}, with using the character table of $G$ in order to compute the structure constants of $G$.
Unfortunately, \cite{Chevie} only contains a partial character table, missing $10$ character families.
However, we can show that these irreducible characters take value $0$ on $T\setminus\{1\}$. In fact, looking at \cite{Lu}, one can verify that the missing characters have the following degrees:
$$\begin{array}{rcl}
\eta_1 & = & t t_2 (q^4-q^2+1) (q^2-\sqrt{2}q+1) (q+1) (q-1) (q^2+1)^2,\\ 
\eta_2 & = & t t_2 (q^4-q^2+1) (q^2+\sqrt{2}q+1) (q+1) (q-1) (q^2+1)^2, \\ 
\eta_3 & = & \frac{\sqrt{2} q}{2} (q+1)(q-1) \eta_1,\\ 
\eta_4 & = & \frac{\sqrt{2} q}{2} (q+1)(q-1)\eta_2,\\ 
\eta_5 & = & (q^2-\sqrt{2} q+1) (q+1) (q-1) \eta_1,\\ 
\eta_6 & = & (q^2+\sqrt{2} q+1) (q+1) (q-1)\eta_2,\\ 
\eta_7 & = & q^4 \eta_1,\\ 
\eta_8 & = & q^4 \eta_2 
  \end{array}$$
where $t_2=q^4-\sqrt{2} q^3+q^2-\sqrt{2} q+1$. Note that $|G:T|=q^{24} t_2 (q^8-1)(q^6+1)(q^4+1)(q^2-1)$ and that $\gcd(|T|,|G:T|)=1$.
Now, let $r$ be a prime dividing $t$. By the above, $r$ does not divide $\frac{|G|}{\eta_k}$ and so all the missing characters are of $r$-defect $0$, and hence vanish at the $r$-singular elements.
So, we can still use \cite{Chevie} for our computations.

Let $c$ be a conjugacy class consisting of regular unipotents ($4$ classes)  and let $c_{3}$ be a class of elements of order $t$ (with the parameter  $\mathtt{aA}=1$).
Computing $\Delta_G(c,c,c_{3})$ we obtain the value $\frac{1}{16} q^{17} t b$, where
$$\begin{array}{rcl}
b & =&  q^{23}-\sqrt{2} q^{22}+q^{19}-\sqrt{2} q^{16}-q^{15}+\sqrt{2} q^{14}+2 q^{13}-2 q^{11}-2\sqrt{2} q^{10}+2 q^9\\
&&+ 3\sqrt{2} q^8 -q^7
+ 6\sqrt{2}q^6+14 q^5+18\sqrt{2} q^4+33 q^3+ 16\sqrt{2} q^2+14q+ 8\sqrt{2}.
  \end{array}
$$
Since $\Delta_G(c,c,c_{3})>0$, we have the desired result. 
\end{proof}

\begin{lemma}\label{Sp4odd}
For every regular unipotent element $u$ of $G=Sp_4(q)$, $q$ odd, there exists $g \in G$ such that $G=\lan u,u^g\ran$.
\end{lemma}

\begin{proof}
As before, we can  compute the structure constants using the character table of $G=Sp_4(q)$ (see \cite{Bl,Pr,Sri}).
We follow the notation of \cite{Sri}. So, let $c$ be any of the two conjugacy classes of the regular unipotent elements $A_{4,1}$, $A_{4,2}$, and let $c_3$ be the conjugacy class of $B_1(1)$, a semisimple element of order $q^2+1$. 
The only irreducible characters taking nonzero values on both classes $c$ and $c_3$ are  
$\chi_1(j)$, $\theta_7$ and $\theta_8$. So,
$$\Delta_G(c,c,c_3)=\frac{q^5-q}{4}\left( q+(-1)^{(q-1)/2}\right),$$
which is positive for all $q\geq 3$.
This implies that, given any regular unipotent element $u$, there exists $g \in G$ such that $uu^g$ has order $q^2+1$.
Note that the Jordan form of $u$ is $J_4$.

If $q=3$, there is no maximal subgroup containing elements of both orders $9$ and $10$. So, the result easily follows.
Then suppose $q>3$. By \cite{MSW} and \cite{BHRD}, the only maximal subgroup containing elements of order $q^2+1$ is $M=H:2$, where $H\cong Sp_2(q^2)$.
So, assume that $X=\lan u, u^g\ran$ is contained in $M$. Since the unipotent elements of $SL_2(q^2)$ have Jordan form $J_2$, we easily get a contradiction, proving that $X=Sp_4(q)$.
\end{proof}

For a fixed $g_3 \in c_3$, denote by $\Delta^*_G(c_1,c_2,c_3)$ the number of distinct triples $(g_1,g_2,g_3)$ such that 
$g_1 \in c_1$, $g_2\in c_2$, $g_1g_2=g_3$, and  $G=\langle g_1,g_2\rangle$. We aim to show that $\Delta^*_G(c_1,c_2,c_3)$ is positive for  certain classes $c_1,c_2,c_3$.
To this end, let $H$ be a maximal subgroup of $G$ containing a fixed element $g_3\in c_3$, and denote by $\Sigma_H(c_1,c_2,c_3)$ the number of distinct pairs 
$(g_1,g_2)\in c_1\times c_2$ such that $g_1g_2=g_3$ and $\langle g_1,g_2\rangle \leq H$. The value of $\Sigma_H(c_1,c_2,c_3)$ can be obtained as the sum of the structure constants $\Delta_H(\tilde c_1, \tilde c_2, \tilde c_3)$ of $H$ for all the $H$-conjugacy classes $\tilde c_1,\tilde c_2,\tilde c_3$ such that $\tilde c_i\subseteq H\cap c_i$. 

The following holds (e.g., see \cite{GM}):

\begin{lemma}
Let $G$ be a finite group and let $H$ a subgroup of $G$ containing a fixed element $x$. Denote by $h(x,H)$ the number of the distinct conjugates of $H$ containing $x$. If  $\gcd(|x|,|N_G(H):H|)=1$, then
$$h(x,H)=\sum_{i=1}^s\frac{|C_G(x)|}{|C_{N_G(H)}(x_i)|},$$
where $x_1,\ldots,x_s$ are representatives of the $N_G(H)$-conjugacy classes fused to the $G$-class of $x$.
\end{lemma}

As a consequence, we obtain a useful lower bound for $\Delta_G^*(c_1,c_2,c_3)$. Namely:
$$\Delta_G^*(c_1,c_2,c_3)\geq \Theta_G(c_1,c_2,c_3),$$ where 
$$ \Theta_G(c_1,c_2,c_3)= \Delta_G(c_1,c_2,c_3) -\sum h(g_3,H)\Sigma_H(c_1,c_2,c_3),$$ 
$g_3$ is a representative of the class $c_3$, and the sum is taken over the representatives $H$ of the $G$-classes of maximal subgroups of $G$ containing elements of all the classes $c_1,c_2,c_3$.  
For groups of small order,  $\Theta_G(c_1,c_2,c_3)$ can be computed using \cite{GAP}. 

\bl{U4(2)}
For every regular unipotent element $u$ of $G=SU_4(2)$, there exists $g \in G$ such that $G=\lan u,u^g\ran$.
\el

\bp
As explained before, the result for $G=SU_4(2)$ can be obtained by using the character table, showing that $\Theta^*_G(c,c,c_3)>0$ for a class $c$ consisting of regular unipotents and a suitable class $c_3$.

The regular unipotent elements of $G$ have Jordan form $J_4$ and belong to the class $4b$. 
We have $\Delta_G(4b,4b,9a)=486$. 
The unique maximal subgroup of $G$ containing elements of both classes $4b$ and $9a$ is $H=3^3:S_4$. Furthermore, 
$h(g_3,H)=1$ and $\Sigma_H(4b,4b,9a)=81$, where $g_3$ is a fixed elements of $9a$.
So, $\Theta_{G}(4b,4b,9a) =  405$, which implies the result.
\enp

\bl{F4(2)}
For every regular unipotent element $u$ of $G=F_4(2)$, there exists $g \in G$ such that $G=\lan u,u^g\ran$.
\el

\bp
The regular unipotent elements of $G=F_4(2)$ have order $16$ and belong to the classes $16a,16b,16c,16d$. 
Using the character table of $G$, we have $\Delta_G(16a,16a,13a)=\Delta_G(16b,16b,13a)=808763850752$ and $\Delta_G(16c,16c,13a)=\Delta_G(16d,16d,13a)=808582807552$.
The list of the maximal subgroups of $G$ were determined in \cite{NW}.
The maximal subgroups of $G$ of order divisible by $13$ do not contain elements from the classes $16a$ or $16b$.
On the other hand, the unique maximal subgroup of $G$ containing elements of classes $13a$ and  $16X$ ($X\in \{c,d\}$) is $H={}^2F_4(2)$.  Furthermore, 
$h(g_3,H)=1$ and $\Sigma_H(16c,16c,13a)=\Sigma_H(16d,16d,13a)=519168$, where $g_3$ is a fixed elements of the class $13a$.
So, $\Theta_{G}(16c,16c,13a) = \Theta_{G}(16d,16d,13a) =  808582288384$, which implies the result.\enp

\begin{rem}\label{Tits}
For $G={}^2F_4(2)'$ the elements of order $16$ are regular unipotent in the simple group of type $F_4$.
We have $\Theta_G(c,c,13a)=\Delta_G(c,c ,13a)=64896$ for $c \in \{16a, 16b,16c, 16d \}$. So, the simple group ${}^2 F_4(2)'$ is generated by two conjugate regular unipotents.
\end{rem}

\begin{lemma}\label{qq3}
Let $G$ be one of the groups $SL_3(q)$ with $q\geq 2$, or $SU_3(q)$ with $q\geq 3$.
For every  regular unipotent element $u$ of $G$,  there exists $g \in G$ such that $G=\lan u, u^g\ran$.
\end{lemma}

\begin{proof}
If $q\leq 5$, we can compute with \cite{GAP} the values of $\Theta_G(c,c,c_3)$ for every conjugacy class $c$ of regular unipotent elements and a suitabile class $c_3$:
$$\begin{array}{c|c|c|c||c|c|c|c}
G & c & c_3 & \Theta_G(c,c,c_3) & G & c & c_3 & \Theta_G(c,c,c_3) \\ \hline
SL_3(2) & 4a & 7a & 7  &  & &  & \\
SL_3(3) & 3b & 13a & 39  & SU_3(3) & 3b & 7a & 56 \\
SL_3(4) & 4a,4b,4c & 7a & 56  & SU_3(4) & 4a & 13a & 247 \\
SL_3(5) & 5b & 31a & 589  & SU_3(5) & 5b,5c,5d & 7a & 7
\end{array}$$
Since $\Theta_G(c,c,c_3)>0$, the group $G$ is generated by two conjugate regular unipotents.

Now, suppose $q\geq 7$ and set $\delta=+1$ if $G=SL_3(q)$ and $\delta=-1$ if $G=SU_3(q)$. Also, set $d=\gcd(q-\delta,3)$.
In the notation of \cite{SF} (see also \cite{G, o13}), the regular unipotent elements belong to the $d$ classes 
$C_3^{(0,\ell)}$, where $0\leq \ell< d$.
So, let $c$ be any of these classes,  and let $c_3$ be the conjugacy
class of $C_8^{(1)}$. Then $c_3$ consists of elements of order
$q^2+\delta q +1$. Using the character table of $G$, one can easily compute that
$\Delta_G(c,c,c_3)=\frac{(q^2+\delta q +1)(q^2-d\delta q-1)}{d^2}$.
Since $\Delta_G(c,c,c_3)>0$, there  exists $g \in G$ such that
the product $uu^g$ has order $q^2+\delta q +1$.

Consider the subgroup $H=\lan u, u^g\ran$. 
Since $H$ contains an element of prime order $p$ (where $p\mid q$) and an element of order $q^2+\delta q +1$, looking at the list of the maximal subgroups of $G$, we obtain that either $H=G$ or $p= 3$ and $H$ is conjugate to a maximal subgroup $A:3$ in $\mathscr{C}_3$ (a unique class of subgroups), where $A$ is a cyclic group of order $q^2+\delta q+1$.
Hence, if $p\neq 3$ we are done.

Finally, assume $p=3$. There are a unique class $c$ of regular unipotent elements and $(q^2+\delta q)/3$ conjugacy classes of elements of order $q^2+\delta q +1$.
By Ito's theorem,  the irreducible characters of $H$ are linear (three characters) or have degree $3$ ($(q^2+\delta q)/3$ characters).
The irreducible characters of degree $3$ have $3$-defect zero, so they take non-zero values only on $A$.
Hence, it is quite easy to compute the structure constants. Let $\tilde c_1,\tilde c_2$ be the two conjugacy classes of $H$ consisting of elements of order $3$ and let $\tilde c_3$ be a conjugacy class consisting of elements of order $q^2+\delta q +1$. For $i,j=1,2$, we have
$\Delta_H(\tilde c_i,\tilde c_j, \tilde c_3)=(q^2+\delta q +1)(1-\delta_{i,j})$.
Given $g_3 \in H$ of order $q^2+\delta q +1$, we have $h(g_3, H)=1$ and then
$\Theta_G(c,c,c_3)\geq (q^2+\delta q +1)(q^2-\delta q-1)-2(q^2+\delta q +1)=(q^2+\delta q +1)(q^2-\delta q-3)>0$.
\end{proof}

\begin{lemma}\label{SU32}
The group $SU_3(2)$ is generated by two regular unipotents. 
\end{lemma}

\begin{proof}
The group $G=SU_3(2)$ admits three conjugacy classes of regular unipotent elements. 
We can compute with \cite{GAP} the character table of $G$ and the structure constant $\Delta_G(4a, 4b, 4c)=10$.
Looking at the maximal subgroups of $G$, only the parabolic subgroups have at least three conjugacy classes of elements of order $4$.
Hence, we can also easily compute $\Theta_G(4a,4b,4c)=8$. This implies that there exist two regular unipotent elements $u,v$ such that $G=\lan u,v\ran$.
\end{proof}

\begin{lemma}\label{2gg}
Let $G=G_2(q)$  with $q\in \{3,4,5\}$.
For every  regular unipotent element $u$ of $G$,  there exists $g \in G$ such that $G=\lan u, u^g\ran$.
\end{lemma}

\begin{proof}
We can compute with \cite{GAP} the values of $\Theta_G(c,c,c_3)$ for every conjugacy class $c$ of regular unipotent elements and a suitable class $c_3$:
$$\begin{array}{c|c|c|c}
G & c & c_3 & \Theta_G(c,c,c_3)  \\ \hline
G_2(3) & 9a & 13a & 7293\\
G_2(3) & 9b,9c & 13a & 7410\\
G_2(4) & 8a & 13a &  245440 \\
G_2(4) & 8b & 13a &  241540 \\
G_2(5) & 25a & 31a &   9373625 
\end{array}$$
As $\Theta_G(c,c,c_3)>0$, the group $G$ is generated by two conjugate regular unipotents.
\end{proof}

\bl{G27}
Let $G=G_2(q)$  with $q\geq 7$ odd.
For every  regular unipotent element $u$ of $G$,  there exists $g \in G$ such that $G=\lan u, u^g\ran$.
\el

\bp 
The maximal subgroups of $G=G_2(q)$ are described in \cite{BHRD}. In particular, as shown in \cite{We}, the only proper maximal 
subgroups of $G$ containing a cyclic torus $T$ of order $t=q^2-q+1$ are of type  $SU_3(q).2$.

Firstly we show that, given a unipotent element $u\in G$, there exists $x \in G$ such that $u u^x$ has order $t$.
This can be achieved by computing the  structure constants of $G$ with \cite{Chevie}.
We follow the notation of \cite{Chevie}, distinguishing two cases.

Suppose $q\equiv 0 \pmod 3$.
There are three conjugacy classes of regular unipotent elements (labeled by $c_{7}$,  $c_8$ and $c_9$, whose centralizers have order $3q^2$).
The class $c_{28}$ (with the parameter \texttt{iI}=1) consists of elements of order $t$.
We have 
$$\Delta_G(c_j,c_j,c_{28})=\left\{ \begin{array}{ll}
\frac{1}{9} q^3 (q^7-q^5-2q^4+8q^3-10q^2+6q-1) & \textrm{ if } j=7,\\[3pt]
\frac{1}{9} q^3 (q^7-q^5-5q^4+2q^3-7q^2-4) & \textrm{ if } j=8,9.
\end{array}\right.$$
Next, suppose that $q$ is odd and such that $q\equiv 1,2 \pmod 3$.
There is only one conjugacy class of regular unipotent elements (labeled by $c_{7}$ whose centralizer has order $q^2$).
The class $c_{31}$ (with the parameter \texttt{iI}=1) consists of elements of order $t$.
We have 
$$\Delta_G(c_7,c_7,c_{31})=\left\{ \begin{array}{ll}
q^3 (q^7-q^5-2q+1) & \textrm{ if }  q\equiv 1 \pmod 3,\\[3pt]
q^5 (q^5-q^3+q^2-1) & \textrm{ if }  q\equiv 2 \pmod 3.
\end{array}\right.$$
Since all these values are positive, we obtain that there exists $x \in G$ such that $u u^x$ has order $t$.
This implies that the subgroup $H=\langle u,u^x\rangle$ is contained in  a maximal subgroup $X=SU_3(q).2$. 

As $q$ is odd, $u$ lies in  $X_{\un}$. By \cite[Theorem 1.1]{McN}, if $q\neq 9$ then $X_{\un}$ is completely reducible.
As $u\in G \leq \Omega_7(q)$ is similar to $J_7$, $X_{\un}$ must be irreducible. However, by \cite[Theorem 1.9]{Sup} $X_{\un}$ has no irreducible representation of degree greater than $3$ containing $J_7$. It follows that $H=G$.
The same conclusion can be obtained also for $q=9$ since the group $SU_3(9)$ does not contain elements of order $|u|=9$.
\enp

\begin{lemma}\label{2s6}
The group $Sp_6(9)$ is generated by two conjugate regular unipotents. 
\end{lemma}

\begin{proof}
Let $G=Sp_6(9)$ be defined by the Gram matrix 
$a=\begin{pmatrix} 0 & \Id_3 \\ -\Id_3 & 0\end{pmatrix}$.
Take $x=\diag(J_3, {}^t J_3^{-1})$ and $y=
\begin{pmatrix}
  1 &  0 &  0 &  0 &  0 &  0\\
  0 &  1 &  0 &  0 &  0 &  0\\
  0 &  0 &  1 &  0 &  0 &  \beta\\
  0 &  0 &  0 &  1 &  0 &  0\\
  0 &  0 &  0 &  0 &  1 &  0\\
  0 &  0 &  0 &  0 &  0 &  1
\end{pmatrix}$,
where  $\beta \in \FF_9$ is such that $\beta^2=\beta+1$.
Then $x,y$ are elements of $G$ and $u=xy$ is regular unipotent (a general construction will be described in Lemma \ref{pa2}).
Let $v=h^{-1}u h$, where $h=(ax)^2a$.
Then, $|uv|=120$, $|uvu|=164$, $|uvu^2|=146$ and $|uv^2u| =728$.
This implies that the order of $H=\langle u,v\rangle$ is divisible by same primes dividing the order of $G$.
By \cite[Corollary 5]{LPS} it follows that $G=H$.
\end{proof}

\begin{lemma}\label{SU43}
The group $SU_4(3)$ is generated by two conjugate regular unipotents.
\end{lemma}

\begin{proof}
Let $G=SU_4(3)$ be defined by the antidiagonal matrix $a=\mathrm{antidiag}(1,1,1,1)$.
Take $u=\begin{pmatrix}
 0  & 1&   0&   0\\
 0  & 0 &  0 & \beta^2\\
\beta^2& \beta^6 &   0  & 0\\
0  & 0 &  1&   1
\end{pmatrix}$
and $x=\diag(\beta^7, \beta^5, \beta, \beta^3)$,
where  $\beta \in \FF_9$ is such that $\beta^2=\beta+1$.
Then $u,x$ are elements of $G$ with $u$ regular unipotent.
Let $v=x^{-1}u x$.  Then, $|uv^6|=28$   and $|u^2v^6| =36$.
Looking at the list of the maximal subgroups of $G$, it follows that $\langle u,v\rangle=G$.
\end{proof}

\begin{lemma}\label{SU53}
The group $SU_5(3)$ is generated by two conjugate regular unipotents. 
\end{lemma}

\begin{proof}
Let $G=SU_5(3)$ be defined by the antidiagonal matrix $a=\mathrm{antidiag}(1,1,1,1,1)$.
Take $u=
\begin{pmatrix}
     1 &      1 &     0 &     0 &     0\\
     0 &      1 &     1 &     1 &     2\\
     0 &      0 &     1 &     2 &     1\\
     0 &      0 &     0 &     1 &      2\\
     0 &      0 &     0 &     0 &      1
\end{pmatrix}$
and $x= 
\begin{pmatrix}
     0&      0&      0&      0  &    1\\
     0&      0&      0&      1&      0\\
     0&      0&      1&      0&      0\\
     0&      1&      0&   \beta^2 &      0\\
     1&      0&      0&      0&  \beta^2
\end{pmatrix}$,
where  $\beta \in \FF_9$ is such that $\beta^2=\beta+1$.
Then $u,x$ are elements of $G$ with $u$ regular unipotent.
Let $v=x^{-1}u x$.  Then, $|uv|=80$, $|uv^3|=61$  and $|uv^5| =28$.
This implies that the order of $H=\langle u,v\rangle$ is divisible by same primes dividing the order of $G$.
By \cite[Corollary 5]{LPS} it follows that $G=H$.
\end{proof}

We combine the results of this section  in a uniform statement as follows.

\begin{propo}\label{pp0}
The following finite quasisimple groups of Lie type are generated by two conjugate regular unipotents:
$$SL_3(q),\; SU_3(q),\; Sp_4(q)\; \text{($q$ odd)},\; Sz(2^{2m+1}),\; {}^2G_2(3^{2m+1}),\; {}^2F_4(2^{2m+1}),\; {}^3D_4(q),$$ $$G_2(q)\; \text{($q$ odd)},\; SU_4(2),\; SU_4(3),\; SU_5(2),\; SU_5(3),\; Sp_6(9),\; G_2(4),\;F_4(2),\; {}^2F_4(2)',\; {}^2E_6(2).$$ 
\end{propo}

Theorem \ref{t22} now follows from Lemma \ref{gt1}, Corollary \ref{co1} and Proposition \ref{pp0}.

\section{Generation by three regular unipotents}

In this section we prove Theorem \ref{qs3}, i.e, that every finite quasisimple group of Lie type is generated by three regular unipotent elements.
In several cases, we show that these three elements can be  chosen to be conjugate, proving Theorem \ref{t33}.

By \cite[Lemma 3.1]{GS} the groups $SL_2(q)$, $q>2$ even, and $\SL_2(9)$ are 
generated by three conjugate regular unipotents. In addition, we have:

\bl{+23} The group $SL_n(q)$, $n>2$, is generated by three  conjugate regular unipotents.\el

\bp The group $SL_3(q)$ is generated by two regular unipotents by Lemma \ref{qq3}. 
Suppose that $n>3$. 
Let $P$ be the stabilizer of a line and $P_1\leq P$ the stabilizer of a vector. 
Then $P_1=O_p(P)\cdot L_1$, where $L_1\cong SL_{n-1}(q)$. 
Then  $L_1$ is generated by two of its regular unipotents. 
As $P=O_p(P)\cdot L$, where $L\cong GL_{n-1}(q)$, these are conjugate in $L$, and hence 
$P_1$ is generated by two conjugate regular unipotents 
$u,u'$ of $G$. 
Let $v\notin P$ be a conjugate of $u$. Then $\lan v,u,u'\ran=G$  by  Lemma \ref{nn6}.
\enp

\bl{ab7}   
Let $G=GL_n(q)$ and $S=SL_n(q)$, $n>1$. Then every \ir $\FF_qG$-module is \ir on $S$. 
\el

\bp  
Suppose the contrary. Then, by Clifford's theorem, $V=V_1+\cdots+V_k$, where   $V_1\ld V_k$ are \ir $\FF_qS$-modules, and $V_i=V_1^{g_i}$ for some $g_i\in G$. 
As $\FF_q$ is a splitting field for $G$ (see \cite[Prop 5.4.4]{KL} or \cite[Lemma 8.5]{ga}), the modules $V_1\ld V_k$ are absolutely irreducible. By Steinberg's theorem, $V_1$ extends to $SL_n(\overline{\FF}_q)$, and we can assume that $V_i=V_1^{g_i}$ as $SL_n(\overline{\FF}_q)$-modules. These are still non-equivalent to each other. As $GL_n(\overline{\FF}_q)=Z\cdot SL_n(\overline{\FF}_q)$, where $Z$ the group of non-zero scalar matrices in $G$, we have $g_i\in Z\cdot SL_n(\overline{\FF}_q)$, and hence the conjugation by $g_i$ yields an inner \au of $SL_n(\overline{\FF}_q)$. Then $V_1,V_i$ are equivalent $SL_n(\overline{\FF}_q)$-modules. This is a contradiction.   
\enp

\bl{a88}\cite[Lemma 2.6]{GS}
Let $G$ be a finite group and $\phi: G\ra GL_m(q)$ be a non-trivial  absolutely \irr of $G$. Suppose that $q$ is minimal in the sense that $\phi$ is not equivalent to a \rep 
into $GL_{m}(q_0)$ for any $q_0\mid q$. Let $V$ be the underlying $\FF_qG$-module. 
Then $V$ is \ir as $\FF_pG$-module (where $p\mid q$ is a prime).\el

\bl{uu6} Let $u\in GL_n(F)=GL(V)$ be a unipotent element. Let $V_1\subset V_2$
be $u$-stable subspaces of $V$.  Suppose that $[u,V_2]\subset V_1$ and $\dim V_2/V_1=d$. Then the Jordan form of $u$ has at least $d$ Jordan blocks, equivalently, $\dim V/[V,u]=\dim C_V(u)\geq d$.   \el

The following lemma is obvious.

\bl{qu9} 
Let $G$ be a quasisimple group of Lie type, $P$ a parabolic subgroup of $G$, and $U$ the unipotent radical of $P$.  Let $V:=Q_1< Q_2\leq U$ be a composition factor of $P$ on $U$. Then $V$ is an elementary abelian $p$-group.\el

In the notation of previous lemma, assume  $G\neq {}^2B_2(q^2), 
{}^2G_2(q^2),{}^2F_4(q^2)$. Since $V$ is an elementary abelian $p$-group,  
$V$ is an \ir $\FF_pL$-module for every Levi subgroup $L$ of $P$. The fact going back to \cite{abs} is that $V$ is completely reducible, in fact, an \ir $\FF_qL$-module in the cases arising below.

\begin{propo}\label{abs}\cite[Theorem 1]{abs}
Suppose that $G$ is a split Chevalley group over a field $\mathbb{K}$, except those arising from a simple algebraic group of type $B_n,C_n,F_4,G_2$ in characteristic $2$ and $G_2$ in characteristic $3$.   Let $P=QL$ be a parabolic subgroup of $G$ with $U$ the unipotent radical of $P$ and $L$ a Levi factor. 
\begin{itemize}
\item[$(1)$] {\rm  \cite[Theorem 1]{abs}} The  quotients of subsequent terms of the central series of $U$ have the structure of completely reducible $\mathbb{K}L$-modules with the simple constituents $M$ being highest weight modules. The highest weights are independent of $\mathbb{K}$ and each simple module is a chief factor of $P$.

\item[$(2)$] {\rm  \cite[Theorem 2(c)]{abs}}. Let $L_{\un}$ be the subgroup of $L$ generated by  the unipotent elements of $L$. If $L_{\un}$ acts on $M$ non-trivially, then $M$ is \ir as an $\FF_qL_{\un}$-module.

\item[$(3)$] {\rm  \cite[Lemma 7(b)]{abs}} If $G$ is non-twisted and    $L_{\un}$ acts on $M$ non-trivially, then $M$ is \ir as $\FF_pL_{\un}$-module.
\end{itemize}
\end{propo}

In the following lemma we exclude the groups $SU_4(3)$ and $ SU_5 (3)$ as these are generated by two conjugate regular unipotents by Lemma \ref{SU43} and  \ref{SU53}. 
 
\bl{unn} 
Let $G=SU_n(q)$, where either $n>5$, or $n=4,5$ and $q\neq 3$ is odd. Then $G$ is generated by three regular unipotents. If $n$ is odd then these can be chosen conjugate to each other. \el

\bp  Let $V$ be the underlying $\FF_{q^2}G$-module,  $W$ be a maximal totally isotropic subspace of $V$ and $m=\dim W$. Let  $P$ be the stabilizer of $W$ in $G$ and  $U$ be the unipotent radical of $P$. Then $P=UL$, where $L$ is a Levi subgroup of $P$. Then  $L_{\un}\cong SL_m(q^2)$. Let $g\in P$ be a regular unipotent and  $X=\lan L_{\un}, g\ran$. By Lemma \ref{qg1}, $g\notin L$, and hence $U_1=X\cap U\neq 1$.

Suppose first that $n$ is even. Then $m=n/2$ and $U$ is an elementary abelian subgroup. By \cite[Lemma 4.6$(1)$]{EZ}, $U$ is \ir as an $\FF_p L_{\un}$-module. This means that $U$ has no proper non-trivial $L_{\un}$-invariant subgroup. 
As $U_1$ is $L_{\un}$-invariant, we have $U_1=U$ and $X=U L_{\un}$. Lemma \ref{gt1} implies that $X=\lan g,h\ran$ for a regular unipotent
element $h\in P$, so $X$ is generated by two regular unipotents. 

Next, let $y$ be a conjugate of $g$ that is not in $P$, and let $Y=\lan g,h,y\ran$. Then
$UL_{\un} <  Y$ and hence $Y$ contains a \syl of $G$ for $p$ dividing $q$. By Lemma~\ref{nn6}, $Y=G$.

Suppose  now that $n=2m+1$ is odd. Then $m=(n-1)/2$.  Let $R=C_P(W)$. Choose a basis  $\{b_1\ld b_n\}$ of $V$ such that  
$(b_{m+1+j},b_i)=\delta_{i,j}$ for $1\leq i,j\leq m$, $(b_{m+1},b_{m+1})=1$ and $(b_{m+1},b_i)=(b_{m+1},b_{m+1+i})=0$.
Then the elements of $R$ are of shape 
$\begin{pmatrix} \Id& {}^tA&B\\0&1&-\overline{A}\\0&0&\Id   \end{pmatrix}$, 
where $B$ is an $(m\times m)$-matrix and $A$ is an  $(1\times m)$-matrix. So, $R$ is a unipotent normal subgroup of $P$, in fact, $R=U$. Let 
$Q$ be a subgroup of $U$ consisting of  the above matrices with $A=0$. Then the action of $P$ on  $U$ by conjugation leaves $Q$ invariant and $U/Q$ is an $\FF_{q_2}L$-module. Note that $L=  \diag\left(h\up,\det\left(h\; {}^t\overline{h}\up\right), {}^t \overline{h}\right)$, where $h$ ranges over $GL_m(q^2)$. 
Then $L\cong GL_m(q^2)$ and $U/Q$ can be viewed as $\FF_{q_2}GL_m(q^2)$-module with natural action of $SL_m(q^2)$. The action of $L\cong GL_m(q^2)$ on $Q$ is as above. 

Observe that $U$ is not abelian. If $U_1\cap Q=1$ then the projection $\overline{U}_1$ of $U_1$ into $Q$ is non-trivial. As $U/Q$ is \ir as $\FF_p SL_m(q^2)$-module, we conclude that $\overline{U}_1=U/Q$. As $[U,U]\neq 1$, we conclude that $U_1\cap Q\neq 1$. As above, this implies $Q\leq U_1$.

 We wish to show that $U_1\cap Q\neq Q$. For this it suffices to observe that $g\notin LQ$. With respect to the above basis, $Q$ consists of matrices
 $\begin{pmatrix} \Id&0&B\\0&1&0\\0&0&\Id   \end{pmatrix}$ and 
 $L=\diag\left(h\up,\det\left(h\; {}^t\overline{h}\up\right), {}^t \overline{h}\right)$.  
\itf $Qb_{m+1}\subseteq \lan b_{m+1}\ran$ and $Lb_{m+1}=  b_{m+1}$. Therefore, $\lan b_{m+1}\ran$ is stable under $LQ$. In addition, if $x\in LQ$ is unipotent then $LQW=W$
implies that $x$ stabilizes a line of $W$. Therefore, $g\in LQ $ implies that $g$
stabilizes two distinct lines on $V$, which is false as $g$ is similar to $J_n$.   

In addition, if $n$ is odd then $L\cong GL_{(n-1)/2}(q^2)$. 
Indeed, in this case $V=W+W_1+W'$, where $W,W'$ are totally isotropic of dimension $(n-1)/2$ and $W_1=(W+W')^\perp$. Then $L$ can be chosen  to be the stabilizer of 
$W$ and $W'$ in $G$. There is a basis with Gram matrix $\begin{pmatrix}0&0&\Id_m\\0&1&0\\
\Id_m&0&0\end{pmatrix}$ with $m=(n-1)/2$, such that $L$ consists of all matrices 
$\diag\left(M\up ,\det\left(M\;{}^t\gamma (M)\up\right),{}^t\gamma(M)\right)$, where $\gamma $ is the Galois automorphism of $\FF_{q^2}/\FF_q$ and $M\in GL_m(q^2)$. \itf the regular unipotent elements of $P/O_p(P)\cong GL_m(q^2)$ are conjugate.
So the above reasoning shows that we can choose the elements in question to be conjugate in $G$.
\enp

The previous argument does not work for $q$ even and $n=4,5$. In Lemmas~\ref{uu5} and \ref{uu4} below we  assume $q>2$, as for $q=2$ we have a stronger conclusion, see  Lemmas~\ref{U5(2)} and \ref{U4(2)}.

\bl{uu5} Let $G=SU_5(q)=SU(V)$, $q>2$ even. Then $G$ is generated by three conjugate regular unipotents.\el

\bp Choose a basis $\{b_1\ld b_5\}$ of $V$ such that the associate  Gram matrix is the antidiagonal matrix $\mathrm{antidiag}(1,1,1,1,1)$.
Let $W=\lan b_1\ran$ be a totally isotropic subspace of $V$ and $P$ the stabilizer of $W$ in $G$. Then $P$ is a parabolic with Levi $L$, 
and $L_{\un}\cong SU_3(q)$.
Let $U=O_p(P)$. Then $U$ is of extraspecial type and the only $L_{\un}$-invariant subgroup of $U$ is $Q=Z(U)$. By Lemma \ref{qq3}, $L_{\un}$ is generated by two 
conjugate regular unipotents. 
\itf $L_{\un} U$ is generated by two conjugate regular unipotents. 

Indeed, let $X$ be a subgroup 
of $P$ generated by two conjugate regular unipotents such that $X/(X\cap U)=L_{\un}$.
We know that $U_1:=X\cap  U\neq 1$. Then  $U_1\cap Q\neq 1$ as $U/Q$ is an \ir 
  $\FF_p L_{\un}$-module and  $[U,U]=Q$. Then $Q <  U_1$. 

Then the upper unitriangular matrices of $P$ forms a Sylow $2$-subgroup $S$ of $G$ and $P$. In addition,
$L_{\un}$  consists of matrices of the form $\diag(1,D,1)$, with $D\in SU_3(q)$, and 
$Q$ consists of matrices $\Id+x$, where $x$ is a $(5\times 5)$-matrix with $(1,5)$-entry in $\FF_q$ and $0$ elsewhere. 
Then $P_{\un}=L_{\un} Q$ consist of matrices of the form $ \begin{pmatrix} 1&0&b \\0&D&0\\ 0&0&1   \end{pmatrix}$ with $b\in \FF_{q}$;
in particular the $(1,2)$-entry of a matrix in  $L_{\un} Q$ equals $0$. In contrast,
the $(1,2)$-entry of any regular unipotent in $S$ is non-zero. 
This is a contradiction. 

\itf that $X=P_{\un}$, so $P_{\un}$ is generated by two conjugate regular unipotents  $u,v$, say.  There is a conjugate  $u'$ of $u$ that is not in $P$. 
Then $\lan X,u'\ran=\lan u,v,u'\ran=G$ by Lemma \ref{nn6}.
\enp

\bl{uu4} 
Let $G=SU_4(q)=SU(V)$, $q>2$ even. Then $G$ is generated by three conjugate regular  unipotents.\el

\bp Choose a basis $\mathcal{B}=\{b_1,b_2,b_3, b_4\}$ of $V$ with Gram matrix   $\Id_4$. Let  $h=\diag(a,a,a,a^{-3})\in G=SU_4(q)$,
where $a$ is a primitive $(q+1)$-root of unity and let $x, y\in G$ be permutational matrices that act on the basis $\mathcal{B}$, respectively, as the $4$-cycle $(b_1,b_2,b_3,b_4)$ and the $3$-cycle $(b_2,b_3,b_4)$. Then $x^4=\Id_4\neq x^2$, so $x$ is similar to $J_4$ and hence $x$ is a regular unipotent.

Let $x_1 = (yh) x (yh)^{-1}$.  Set $H=\lan x, x_1\ran$, so $H$ is generated by two conjugate regular unipotents. Since $x^2 ((x x_1)^3) x^{-2}=h^{-4}$ and $\gcd(q+1,4)=1$, we observe that $H$ contains $h$ and also $x^i h x^{-i}$ for $i=2,3,4$. So $D<  H$, where $D$ is the group of diagonal matrices in $G$ under the basis $\mathcal{B}$. Note that $|D|=(q+1)^3$.

In particular, $H$ contains the element $d=\diag( a^{-8}, a^4,  1,  a^4)\in D$ and then it contains $d (x_1 x)^2= y$. Note that $\lan x,y\ran \cong S_4$.  It follows that $H/D\cong S_4$ and $H=N_G(D)$. By  \cite[Table 8.10]{BHRD}, $H$ is a maximal subgroup of $G$. Therefore, $G$ is generated by $x,x_1$ and any conjugate $x_2$
of $x$ such that $x_2\notin H$.
\enp

\bl{oo1}  Let $G=\Om^+_{2n}(q)=\Om(V)$, $n>3$, and let $P$ be the stabilizer of a totally singular subspace of 
dimension $n$. Then $P_{\un}$ is generated by two regular unipotents and $G$ is generated by three regular unipotents. 
If $q$ is even or $n$ is odd then these can be chosen conjugate.\el

\bp Let $U$ be the unipotent radical of $P$, and $L$ a Levi. Then $L$ is isomorphic to
a subgroup of index $\gcd(2,q-1)$ of $GL_n(q)$. \itf the regular unipotent elements 
of $L$ partition in at most $2$ conjugacy classes; if $q$ is even then they are conjugate in $L$. In $n$ is odd then $|GL_n(q):Z(GL_n(q))\cdot SL_n(q)|=\gcd(n,q-1)$ is odd. Therefore, if we view $L$ as a subgroup of $GL_n(q)$, we observe that
$|GL_n(q): L\cdot Z(GL_n(q))|$ is odd, whence $GL_n(q)=L\cdot Z(GL_n(q))$. So  the regular unipotent elements  
of $L$ are conjugate in $L$.

Under a certain basis of $V$ we can write  $U=\begin{pmatrix}\Id_n&M\\ 0&\Id_n\end{pmatrix}$, where $M$ is the set of skew-symmetric matrices if $q$ odd, and symmetric matrices with zero diagonal if $q$
is even. Note that $P_{\un}=L_{\un} U$ and $L_{\un}\cong SL_n(q)$. We know that $M$ is an \ir $\FF_p L_{\un}$-module (e.g., see \cite[Lemma 4.6]{EZ}). As in Lemma \ref{unn}, we observe that there are regular unipotent elements $g,h\in P_{\un}$ that generate a subgroup $X$ such that $X/(X\cap U)\cong L_{\un}$. In addition, $X\cap U\neq 1$. Therefore, $U < X$, and hence $X=P_{\un}$. Now, choose  a conjugate $g'$ of $g$ with $g'\not \in P$. 
As $P_{\un}$ contains a \syl of $G$, we have $\lan P_{\un},g'\ran=G$ by Lemma \ref{nn6}.
\enp

\bl{oo2} 
Let $G= \Om_{2n+1}(q)$ with $n\geq 3$ and $q$ odd, or $\Om^-_{2n+2}(q)$ with $n\geq 3$. 
Let $V$ be the underlying space of $G$ and let $P$ be the stabilizer of a maximal totally singular subspace. Then 
$P_{\un}$ is generated by two conjugate  regular unipotents and $G$ is generated by three conjugate  regular unipotents.
\el

\bp 
Let $W$ be a maximal totally singular subspace of 
$V$. Then $\dim W=n$ and $V=W +V_1+W'$, where $V_1$ is complement of $W$ in $W^\perp$ and $W'$ is a complement of $W^\perp
$ in $V$ and is totally singular. Let $U$ be the unipotent radical of $P$. Under a certain basis of $V$ the group  $Z(U)$ has form $\begin{pmatrix}\Id_n&0&M\\ 0&\Id&0\\ 0&0&\Id_n\end{pmatrix}$, where $M$ is the set of 
skew-symmetric matrices if $q$ odd, and symmetric matrices with zero diagonal if $q$
is even (see, for instance,  \cite[Lemma 4.14]{EZ}). Let $L=\{h\in P: hV_1=V_1,hW'=W'\}$ be the stabilizer of the above decomposition.
Then $L$ consists of matrices
$\diag(x,y,{}^tx\up)$, where $x\in GL_n(q)$ and $y\in SO(V_1)$, 
see \cite[Lemma 4.8]{EZ} for details. As $\dim V_1\leq 2$ and $SO(V_1)$ is abelian, we have $L_{\un}\cong SL_n(q)$.

Note that $Z(U)$ is an elementary abelian $p$-group for a prime $p$ dividing $q$,
so this can be seen as an $\FF_pL_{\un}$-module. By \cite[Lemma 4.6]{EZ}, this is \ir
and hence has no proper $L_{\un}$-invariant subgroup. As in Lemma \ref{unn}, we observe that there are two regular unipotent elements $g,h\in P_{\un}$ 
that generate a subgroup $X$ such that $X/(X\cap U)\cong L_{\un}$. If   $X\cap Z(U)\neq 1$ then $Z(U) < X$. 
Moreover, here we can assume that $g,h$ are conjugate in $G$.
Indeed, as $\overline{g}:=g\pmod U$ is contained in a subgroup isomorphic to $GL_n(q)$ and $\overline{g}$ is a regular unipotent in $GL_n(q)$, where regular unipotent elements are conjugate,  we can choose $h$ to be a conjugate of $g$ in $G$.

Note that $U$ is non-abelian and $U'= Z(U)$, see  \cite[Lemma 4.18]{EZ}. This implies $X\cap Z(U)\neq 1$. Indeed, otherwise $1\neq X \cap U\not\leq Z (U)$, so $(X \cap U)'\neq 1$, a contradiction.

So, we are left with ruling out the case where $X \cap U=Z(U)$. Suppose that
$X=L_{\un} Z(U)$. The matrices of $X$ are 
$$\begin{pmatrix}A&0&0\\ 0&\Id_k&0\\ 0&0&{}^tA\up\end{pmatrix}\cdot \begin{pmatrix}\Id_n&0&M\\ 0&\Id_k&0\\ 0&0&\Id_n\end{pmatrix}=\begin{pmatrix}A&0&AM\\ 0&\Id&0\\ 0&0&{}^tA\up\end{pmatrix},$$
where $A\in GL_n(q)$ and   $k=1$ if $\dim V$ is odd and $k=2$ otherwise. If $y\in L_{\un} Z(U)$ is unipotent then 
$y$ fixes the vectors $v\in V_1$ and some vector in $W$. So, the fixed point space $V^y$ of $y$ has dimension  $\dim V^y\geq k+1$.

(a) If $k=1$ then $G=\Om_{2n+1}(q)$, $q$ odd, and the Jordan form of a regular unipotent of $G$ is $J_{2n+1}$ so $\dim V^y=1$. This is a contradiction.

(b) If $k=2$ then $G=\Om^-_{2n+2}(q)$ and the Jordan form of a regular unipotent of $G$
is $\diag(1,J_{2n+1})$ if $q$ is odd and $\diag(J_2,J_{2n})$ if $q$ is even 
(see \cite[Proposition 3.5]{LSb} and \cite[Theorem 3.1]{GLOB}). In both cases, we have $\dim V^y=2$. So, this is again a contradiction. 

It follows that  $P_{\un}$ is generated by two conjugate regular unipotents  $u,v$, say.  
Now, choose  a conjugate $u'$ of $u$ with $u'\not \in P$. 
As $P_{\un}$ contains a \syl of $G$, we have $\lan u,v,u'\ran=G$ by Lemma \ref{nn6}.
\enp

\bl{tt3} 
Let $G=Sp_{2n}(q)$,   $q$ even, and let $g\in G$ be a $2$-element. Suppose that $g$ is orthogonally indecomposable. Then $\Jor(g)=J_{2n}$ or $\diag(J_n,J_n)$. In addition,  if $|g|=2^k$  then $2^k$ is the minimal $2$-power such that $2n\leq 2^k < 4n$ or  $n\leq 2^k < 2n$, respectively.
\el

\bp 
This is well known, see for instance \cite[Section 6.1 or Theorem 7.3]{LSb}. Note that there are several conjugacy classes of each type of orthogonally indecomposable elements in $G$. The additional claim is straightforward, or see  \cite[Lemma 4.6]{DZ}.
\enp

\bl{op2}\cite[Lemma 3.12]{z24}
Let $G=O^+_{2n}(2)=SO(W)$, $n>2$. Suppose that $g\in G$ is indecomposable on $W$. Then $G$ is generated by  three conjugates of $g$.
\el

Note that $J^2_{2n}$ is not a regular unipotent in $O^+_{2n}(2)$.

\bl{pa2} 
Let $G=Sp_{2n}(q)=Sp(W)$, where $n>2$, 
and let $P$ be a parabolic with Levi $GL_n(q)$.  Then $P$ contains two regular unipotent elements $u,v$ conjugate in $P$ such that $\lan u,v \ran=P_{\un}$. 
\el

\bp 
Let $G$ be defined by the Gram matrix $\begin{pmatrix} 0&\Id_n\\ -\Id_n&0\end{pmatrix}$.
Let $b_1\ld b_{2n}$ be the basis in which this matrix is written. 
Let $W_1=\lan b_1\ld b_n\ran$  and $W_2=\lan b_{n+1}\ld b_{2n}\ran$. Then $W_1,W_2$ are totally isotropic subspaces. Let $P$ be the stabilizer of $W_1$ in $G$. Then $P$ is a maximal parabolic subgroup. Let $U$ be the unipotent radical of $P$. 
Then $U=\begin{pmatrix} \Id_n&M\\ 0&\Id_n\end{pmatrix}$, where $M$ is the set of symmetric matrices. The subgroup $U$ has a unique proper normal subgroup $Q$ of $P$, which has the same form $Q=\begin{pmatrix}\Id_n & \widetilde M\\ 0 & \Id_n\end{pmatrix}$, where $\widetilde M$ is the set of symmetric matrices with zero diagonal (see for instance \cite[Lemma 4.6]{EZ}). The set of the  matrices  
$\begin{pmatrix} A&0\\ 0& {}^t A\up \end{pmatrix} $, $A\in  GL_n(q)$, is a Levi of $P$. Let $x=\begin{pmatrix}J_n&0\\ 0& {}^t J_n\up \end{pmatrix}$, where $J_n$ is the upper triangular Jordan block. 
Then $x-\Id_n$ is nilpotent and $(x-\Id_n)^ib_n=b_{n-i}$, $(x-\Id_n)^ib_{n+1}\in \lan 
b_{n+1+i},\ldots,b_{2n} \ran$ for all $i=1\ld n-1$. 
Let $y=\begin{pmatrix}\Id_n& C\\ 0&\Id_n\end{pmatrix}$, where $C$
is an $(n\times n)$-matrix with $1$ at the $(n,n )$-position and $0$ elsewhere (so $y=\Id_{2n}+C_1$
where $C_1$ has $1$ at the $(n,2n )$-position and $0$ elsewhere). In particular, $yb_i=b_i$ for $i=1\ld 2n-1$ and $yb_{2n}=b_{2n}+b_n$. 
Set $u=xy$. 

We claim that $u$ is regular unipotent. Let $w\in W$ and $uw=w=\sum a_ib_i$, ($a_i\in\FF_q$). Let 
$\bar w$ the projection of $w$ to 
$W/W_1$. As $u$ acts on $W/W_1$ as ${}^t(J_n\up)$, it follows that the fixed point subspace of $u$ on $W/W_1$ is $\FF_q b_{2n}$. So, $a_{n+1}=\ldots =a_{2n-1}=0$. As $ub_{2n}=x(b_{2n}+b_n)=
b_{2n}+b_n+b_{n-1}$, we conclude that  $a_{2n}=0$, and hence $w\in W_1$. 
As  $u$ acts on $W/W_1$ as $J_n$, it follows that $a_2=\ldots=a_n=0$. So, the fixed point space of $u$ on $W$ is one-dimensional, and the claim follows.  

As $P/U \cong GL_n(q)$, the regular unipotents are conjugate in $P/U$. So, by Lemma \ref{gt1},  $P_{\un}/U$
is generated by two conjugate regular unipotents. Therefore,  $L_{\un} < \lan u,v\ran$ for some $P$-conjugate $v$ of $u$. 

We claim that $\lan L_{\un},u\ran =P_{\un}$.
Suppose first that $q$ is odd. Then $M$ is an \ir $\FF_qL$-module (see, for instance, \cite[Lemma 4.6]{EZ}), which remains \ir over $L_{\un}$ by Lemma \ref{a88}. Moreover,  by Lemma \ref{ab7}, $M$ is \ir as  $\FF_p L_{\un}$-module, and hence $U$ has no proper non-trivial $L_{\un}$-invariant subgroup. So $U <  \lan L_{\un},u\ran $, whence the result. (A more conceptual way could be with using \cite[Prop 5.1.3]{Car}.)

Suppose that $q$ is even. As $x\in L_{\un}$, we have $y\in U\cap \lan L_{\un},u\ran$, where   $y\notin Q$ by the choice of $y$.  So the projection of $y$ into $U/Q$ is non-trivial.
As $U/Q\cong M/\widetilde M$ is an \ir $\FF_qL$-module, 
$U/Q$ is irreducible also  as $\FF_qL_{\un}$-module  (Lemma \ref{a88})
and as  $\FF_p L_{\un}$-module (Lemma \ref{ab7}). 
So, we conclude that $U/(U\cap \lan L,u\ran)=U/Q$.
\itf $Q\cap \lan L,u\ran\neq 1$, since $M$ is an indecomposable $\FF_qL$-module   \cite[Lemma 4.6]{EZ}. Again, $\widetilde M$
is irreducible, so $Q\cap \lan L,u\ran=Q$. This implies the result.  
\enp

\bl{spp2} 
The group $Sp_{2n}(q)$, $n>2$, is generated by three conjugate regular unipotent elements.
\el

\bp 
Let $G,P,u,v$ be as in Lemma \ref{pa2}. Then $P_{\un}$ is generated by two conjugate regular unipotents. Let $u'$ be a conjugate of $u$ such that $u'\notin P$.
By Lemma~\ref{nn6}, we have $\lan u',P_{\un}\ran=G$.\enp

\bl{s44}
Let  $G=Sp_4(q)$, $q$ even. Then $G$ is generated by three conjugate regular unipotents.\el 

\bp  If $q=2$, the result follows from Lemma \ref{ss42}. So, we may assume $q>2$. 
Let $x\in G$ be of order $q^2+1$.
Note that $N_G(\lan x\ran )$ contains an element, $a$ say, of order $4$ (see, for instance,
\cite[Table 8.14]{BHRD}, where a group $C_{q^2+1}:4$ occurs, we do not need to decide on its maximality).
The group $\lan x,a\ran $ is of order $4(q^2+1)$ and $\lan x,a^2\ran$ is dihedral of order $2|x|$, where $|x|=q^2+1$ is odd. \itf $\lan x,a\ran $ contains a conjugate $b$ of $a$ such that $\lan a,b \ran = \lan x,a \ran$; in particular, $\lan x,a\ran $ is generated by two conjugate elements of order $4$.  
Observe that the subgroups of order $q^2+1$ are conjugate in $G$.

By Lemma \ref{ss42}, the group $Sp_4(2)$ can be generated by two conjugate elements with Jordan form $J_4$.
We conclude that there are two conjugates $b,c$ of $a$ such that the group $X=\lan a,b,c\ran$ contains an elements $x$ of order $q^2+1$ and 
a subgroup $Sp_4(2)$. 
By inspection of the maximal subgroups of $G$, listed  in \cite{BHRD}, the only maximal subgroups containing elements of order $q^2+1$ are
$Sp_2(q^2):2$ and $O_4^-(q)$. However, these subgroups do not contain $Sp_4(2)$: in fact, this would imply the existence
of a projective representation of $Sp_4(2)'\cong Alt(6)$ of degree $2$ over a field of characteristic $2$, in contradiction with \cite[Proposition 5.3.7]{KL}.
\enp

\bl{G2even} 
If $q>4$ is even, then the group $G_2(q)$  is generated by three conjugate regular unipotents.
\el

\bp As recalled in the proof of Lemma \ref{G27}, the only proper maximal 
subgroups of $G$ containing a cyclic torus $T$ of order $t=q^2-q+1$ are of type  $SU_3(q).2$.
So it suffices to show that, given a unipotent element $u\in G$, there exist $x,y \in G$ such that $u u^x$ has order $t$ and  $u u^y$ has order $t_2=q^2+q+1$.
This can be achieved by computing the  structure constants of $G$ with \cite{Chevie}.

There are two conjugacy classes of regular unipotent elements (labeled by $c_{7}$ and $c_8$, whose centralizers have order $2q^2$).
The classes $c_{26}$  and $c_{27}$ (with the parameter \texttt{iI}=1) consist of elements of respective order $t_2$ and $t$. For all $j=7,8$, we have 
$$\Delta_G(c_j,c_j,c_{26})=\left\{ \begin{array}{ll}
\frac{1}{4}q^3(q^7-q^5+5q^3+7q^2+5q+1) & \textrm{ if } q\equiv 1 \pmod 3 \\[3pt] \frac{1}{4}q^4(q^6-q^4+q^3+5q^2+6q+3) & \textrm{ if } q\equiv 2 \pmod 3 
\end{array}\right.$$
and 
$$\Delta_G(c_j,c_j,c_{27})=\left\{ \begin{array}{ll}
\frac{1}{4} q^4 (q^6-q^4-q^3+5q^2-6q+3)  & \textrm{ if } q\equiv 1 \pmod 3 \\[3pt]
\frac{1}{4} q^3 (q^7-q^5+5q^3-7q^2+5q-1) & \textrm{ if } q\equiv 2 \pmod 3. 
\end{array}\right.$$ Since all these values are positive, we are done.\enp

\bl{f44} 
The group $F_4(q)$, $q$ odd, is generated by three conjugate regular unipotents.
\el

\bp Let $P=P_2$ be a parabolic subgroup of $G=F_4(q)$ generated by the root
subgroups $x_{\al}(t)$ with $\al\in \{- \al_1,-\al_3,-\al_4,\beta\}$,
where $\beta$ ranges over all positive roots and $t\in \FF_q$.
Set $L_{\un}=\lan x_{\al}(t): \al\in \{\pm \al_1,\pm\al_3,\pm\al_4\}\ran$.
Then $L_{\un}$ is a   subgroup of a Levi subgroup of $P$, 
in fact, $L_{\un}\cong SL_2(q)\times SL_3(q)$,  see \cite[Table 8]{C23}.
Set $u=x_{\al_1}(1)x_{\al_2}(1)x_{\al_3}(1)x_{\al_4}(1)$. Then $u\in P$
is a regular unipotent,  see \cite[Lemma~2.1]{T}.

The information on the unipotent radicals of maximal parabolics of $F_4(q)$
for $q$ odd is available in \cite[Table 2]{Ko14b}. In this case, $U=O_p(P)$
is described as follows.  
The lower central series of $U$ has three terms: $1< Q_1< Q_2 < U$,  and these are unique $L_{\un}$-invariant subgroups of $U$. In addition, the consecutive factors are of size $|Q_1|=q^2$, $|Q_2/Q_1|=q^6$ and $|U/Q_2|=q^{12}$.

Let $X=\lan L_{\un},u\ran$. By \cite[Corollary 5.2.3]{Car} we have
$[u,x_{\al_1}(1)]=x_{\al_1+\al_2}(t).y$, where $y\in Y$,
the group generated by $x_\be$ with $\be\succ \al_1+\al_2$, as
$[x_{\al_3}(1),x_{\al_1}(1)]=[x_{\al_4}(1),x_{\al_1}(1)]=1$.
\itf $[u,x_{\al_1}(1)]\in U\setminus Q_2$. Therefore, $(X\cap U)/(X\cap
Q_2)\neq 1$ is a non-trivial  normal subgroup of $L_{\un} U/Q_2$. This implies
$(X\cap U)/(X\cap Q_2)=U/Q_2$.  (Note that $\al_1+\al_2$ is a weight for $L_{\un}$
on $U/Q_2$, see \cite[Table 2]{Ko14b}, strictly speaking, this is the weight of 
the corresponding   module for algebraic group of type $A_1\times A_2$).

Furthermore, $U/Q_2$ has weights $\al_2$ and $\al_1+\al_2+\al_3$ \cite[Table 2]{Ko14b} and hence there are elements
$y_1,y_2\in U$ whose projection into $U/Q_2$ are $x_{\al_2}(1)$ and
$x_{\al_1+\al_2+2\al_3}(1)$. Then $[y_1,y_2]\in Q_2$, and the projection
of $[y_1,y_2]$ into
$Q_2/Q_1$ is $x_{\al_1+2\al_2+2\al_3}(t)$ for $t\in \FF_q$. As above it follows that
$Q_2/Q_1 <  X/Q_1$.

In addition, there is $y_3\in Q_2$ whose projection into 
$Q_2/Q_1$ is  $x_{\al_1+2\al_2+2\al_3+2\al_4}(1)$ \cite[Table 2]{Ko14b}.
Then $[y_2,y_3]=x_{2\al_1+3\al_2+4\al_3+2\al_4}(b)\in Q_1$. It follows
that $Q_1<  X$.

We conclude that $U<  X$, and hence $X=L_{\un}U$. 
As $L_{\un}$ is generated by two  regular unipotents that can be chosen conjugate if $q\neq 9$ (see Remark \ref{SL29} and Lemma \ref{qq3}), it follows that so is $L_{\un} U$.
Choose a conjugate $u'$ of $u$ such that $u'\notin P$. Then $\lan u',X\ran$
contains a \syl of $G$ for $p$ dividing $q$. Using Lemma \ref{nn6}  we
conclude that $\lan u',X\ran=G$.

Let $q=9$. We choose for $P$ a maximal parabolic subgroup whose Levi satisfies
$L_{\un}\cong Sp_6(9)$. By Lemma \ref{2s6}, $L_{\un}$ is generated by 2 conjugate regular unipotents. By \cite[Proposition 4.5]{CKS}, $U=O_3(P)$ is a group of extraspecial type of order $q^{15}$, $Z(U)=U'$,  and $L_{\un}$ acts irreducibly on $U/Z(U)$ and trivially on $Z(U)$. In fact, 
$U/Z(U)$ is an \ir $\FF_3L_{\un}$-module, so $Z(U)$ is a maximal proper normal subgroup of $L_{\un} U$ contained in $U$. Let $U_1=\lan u,v\ran$, where $u,v$ are conjugate regular unipotents whose projections into $P/U$ generate $L_{\un}$. As in the proof of Lemma \ref{ee6} we observe that  $U_1=U$. Let $v'$ be a conjugate of $v$ such that $v'\notin P$. Then $\lan u,v,v'\ran=G$ by Lemma \ref{nn6}. 
\enp

Note that the above method does not work for the case with $G=F_4(q)$, $q$ even.  

\bl{1w1}  \cite[Propositions 4.4 and 4.6]{CKS}  
Let $G=E_6(q),{}^2E_6(q)$ and let $Q$ be the root subgroup of $G$ corresponding to the maximal root. Let $P=N_G(Q)$ and $L=P/U$, where $U$  is  the unipotent radical of $P$. Then $|U|=q^{21}$ and $L\cong SL_6(q),\SU_6(q)$, respectively. 
Moreover, the conjugation action of  $P$ on $U$ defines on  $U/Q$ a structure of   an $\FF_qL$-module isomorphic  to  the  third  exterior
power of the natural $SL_6(q)$-, $\SU_6(q)$-module, respectively,  and the restriction of $U/Q$ to $L$ is an  \ir module of dimension $20$.  \el

\bl{1w2}  Let $G=E_n(q)$, $n=7,8$, and let $P$ be the maximal parabolic subgroup of $G$ corresponding to the maximal root $\al_2$. 
Let $U$  be  the unipotent radical of $P$ and $L$ a Levi. Then $L_{\un}\cong SL_{n}(q)$.
The $P$-composition series of $U$ is unique of length $2,3$ for $n=7,8$, respectively, the factors  are elementary abelian groups. The conjugation action of  $P$ on  $U$ turns the factors to an \ir $\FF_q L_{\un}$-modules of dimensions $7,35$ and $8,28,56$, respectively. These are \ir as $\FF_p L_{\un}$-modules, and hence these are chief factors of $P_{\un}=L_{\un}U$. In addition, the $P_{\un}$-composition series of $U$ is unique.

Let $Q_1<U$, respectively, $Q_1< Q_2<U$ be the non-trivial proper normal subgroups of $P$. Then $Q_1$ is abelian, $[U,Q_2]= Q_1$, and $U'=Q_1,Q_2$, respectively.   In addition, $U$ and $Q_2$ are non-abelian. \el

\bp Let $p\mid q$ be a prime.  By Proposition \ref{abs}(1) and \cite{Ko14b},   the factors of the lower central series of $U$ are non-trivial \ir $\FF_qP$-modules, in fact, \ir  $\FF_qL$-modules. In \cite{Ko14b}  the author computes their dimensions as mentioned.

By Lemma \ref{ab7} or Proposition \ref{abs}(2), these are \ir $\FF_qL_{\un}$-modules. 
We can use Lemma \ref{a88} to conclude that these are \ir $\FF_pL_{\un}$-modules. For this we need to show that $\FF_q$ is the minimal realization field of their realization.   This is the case by \cite[Proposition 5.4.6]{KL}.

\itf  the lower central series  terms are the only normal subgroups of $P$ that lie in $U$,
so there is no normal subgroup $N$ of $P$ such that $[Q_i,U]<N<Q_i$. \enp
 
\bl{ee6} 
Let $G=E_n(q)$, $n\in\{6,7,8\}$, and let $P$ be a maximal parabolic corresponding to a simple root $\al_2$ in notation of {\rm \cite{Bo}}. Then $P_{\un}$ is generated by two regular unipotent elements of $G$. Consequently, $G_{\un}$ is generated by three regular unipotents.
\el

\bp 
Let $L$ be a Levi subgroup of $P$. By general theory, $L_{\un}\cong SL_{n}(q)$. Let $U$
be the unipotent radical of $P$. Let $g,h$ be two regular unipotents in $P$ that generate $L_{\un}$ modulo $U$ and let $X=\lan g,h\ran$. Observe that $X\cap U\neq 1$. Indeed, 
$P=UL$ is a semidirect product and $g\in UL_{\un}=P_{\un}$ as $L/L_{\un}$ is a $p'$-group.
As $g\notin L_{\un}$ (Lemma \ref{qg1}), it follows that $U_1:=X\cap U\supset \lan L_{\un},g\ran\cap U\neq 1$.
Note that $U_1$ is an $L_{\un}$-invariant subgroup of $U$. If $U_1=U$, we are done, so we assume $U_1<U$. If $U_1<Q_1$ then $U_1 $ is normal in $L_{\un}U$. 
As $L_{\un}$ acts irreducibly on $Q_1$, it follows that $U_1=Q_1$, unless the action of 
$L_{\un}$ is trivial which happens only when $n=6$. 

Suppose first that $G=E_6(q)$. By \cite[Propositions 4.4]{CKS},    $U$ is of extraspecial type and $Z(U)$ is a root subgroup $U_\al$, say.  Note that $U_1$ is not a subgroup of $Z(U)$. Indeed, 
as $[L_{\un},Z(U)]=1$   in fact, $U_1< \lan L_{\un},U_\al,U_{-\al}\ran $, which is a semisimple group. This contradicts \cite[Theorem 1.4]{tz13}. So $U_1$ is not a subgroup of $  Z(U)$, and hence the projection of 
$U_1$ into $U/Z(U)$ is non-trivial. By Lemma \ref{1w1}, $U/Z(U)$ is \ir as $\FF_pL_{\un}$-module, so $U/Z(U)=U_1/(Z(U)\cap U_1)$. As $[U,U]=Z(U)$, it follows that $[U_1,U_1]=Z(U)$, and hence $U_1=U$. This proves the result for $E_6(q)$.
    
Let $n>6$. Then $Q_1$ is a non-trivial \ir $\FF_qL_{\un}$-module as well as  $\FF_pL_{\un}$-module by Lemma \ref{1w2}.  Therefore, $Q_1$ has no $L_{\un}$-invariant  proper non-trivial subgroup, so either $U_1\cap  Q_1=1$ or $Q_1\leq  U_1$. 
 
Suppose that  $U_1\neq Q_1$, and let $\overline{U_1}$
be the projection of  $U_1$ into $U/U'$. As above, $U/U'$ has no $L_{\un}$-invariant  proper non-trivial subgroup. If $\overline{U_1}\neq 1$ then  either $\overline{U_1}=U/U'$ or  $\overline{U_1}<U'$. 
In the former case, $U_1$ is non-abelian and hence $1\neq U_1'\leq Q_1$. So $U_1\cap Q_1\neq 1$, and hence $Q_1\leq U_1$ as $Q_1$ has no proper non-trivial $L_{\un}$-invariant subgroup.  So, $U_1=U$. 
Finally, we conclude that $U_1=	Q_1$ or $Q_2$. 

Suppose that $U_1=Q_1$. Then $[u,Q_2] < Q_1$. By Lemma \ref{1w2},   $\dim_{\FF_q}(Q_2/Q_1)=35,28$
for $n=7,8$, respectively. On the other hand, by \cite[Tables 8,9]{Lw}, $u$ has at most $8$ blocks on the adjoint module both of $n=7,8$. 
This contradicts Lemma \ref{uu6}.  
Similarly, if $G=E_8(q)$ and $U_1=Q_2$ then $\dim U/Q_2=56$ by Lemma \ref{1w2}, and we conclude as above. \enp

\bl{2e6} 
The group ${}^2E_6(q)$, $q$ odd, is generated by three conjugate regular unipotents.
\el

\bp Note that the Dynkin diagram of $G={}^2E_6(q)$ as a group with BN-pair is of type $F_4$. 
This is obtained by gluing the nodes $(1,6),(3,5)$ at the Dynkin diagram of the root system of type $E_6$ and the nodes $(2),(4)$ remain unchanged. We use \cite{Ko14},
where the nodes of the obtained $F_4$ diagram  are denoted by $(1+6),(3+5),(2),(4)$. Let $P$
be the  parabolic subgroup of $G$ corresponding to the node $(4)$ at the obtained $F_4$ diagram. Then the Levi $L$ is described as $(SL_2(q)\times PSL_3(q^2)).(q-1)$, and $L_{\un}\cong  SL_2(q)\times PSL_3(q^2)$. The root subgroup of the multiple $SL_2(q)$ is $x_{\al_2}$ of root system of~$G$.   

Let $U$ be the unipotent radical of $P$. Then $|U|=q^{29}$ by \cite{Ko14}, and then $L_{\un}$ is as above by \cite[Table 10]{C23}. (This is to make it sure for our choice 
of parabolic, another choice leads to $|U|=q^{31}$.) 
The lower (and upper) central series of $P$ on $U$ has three non-trivial terms 
$Q_1 <  Q_2  <  U$, say, whose factors  can be viewed as \ir $\FF_qL$-modules of dimensions $2$, $9$ and $18$, respectively.  
By \cite[Table 1]{Ko14},  there are exactly $3$ chief factors of $P$ on $U$, 
so these are the factors of the lower (and upper) central series of $U$.

Let $u\in P$ be a regular unipotent. Then $u$ is a product of positive root elements.
Let $u= x_{\al_1}(1)x_{\al_6}(1)x_{\al_2}(1)x_{\al_4}(1)x_{\al_3}(1)x_{\al_5}(1)$. 
As $[x_{\al_1}(1), x_{\al_6}(1)]=1$ and $[x_{\al_3}(1), x_{\al_5}(1)]=1$,  
we have   $u\in {}^2E_6(q)$ by the above, and $u\in P$. 
As  $x_{\al_1}(1)x_{\al_6}(1)$, $x_{\al_3}(1)x_{\al_5}(1)$, $x_{\al_2}(1) \in L_{\un}$,  it follows that $\lan L_{\un}, u\ran$ contains an element in $U\setminus Q_2$, in fact, this is of the form $x_{\al_4}(1) \cdot y$, where $y$ is a product of root elements $x_\be$ for some roots $\be \succ \al_4$. (That is, $\be -\al_4$ is a positive linear combination of simple roots.)

Let $X=\lan L_{\un},u\ran$. By \cite{Ko14}, there are no intermediate $L_{\un}$-invariant subgroup $Y$, $Q_2<Y<U$. Therefore, $(X\cap U)/(X\cap Q_2)=U/Q_2$. As $U>Q_2>Q_1>1$
is a lower central series of $U$, it follows that $U' <  Q_2$ and $U'$ is not a subgroup of $Q_1$. As above, the is no intermediate $L_{\un}$-invariant subgroup $Y$, $Q_1<Y<Q_2$, and hence $X/(X\cap Q_1)$ contains $Q_2/Q_1$ and hence $U/Q_1$.
Finally, as $|Q_1|=q^2$, in the additive notation, is not a $1$-dimensional $\FF_qL$-module, this is an \ir $\FF_pL_{\un}$-module for $p\mid q$  (see also Lemma \ref{a88}). The latter means that 
there is no intermediate $L_{\un}$-invariant subgroup $Y$, $1<Y<Q_1$. \itf $X$
contains $U$. 
  
If $q\neq 9$ then $L_{\un}$ is generated by two conjugate regular unipotent elements. 
Therefore, $P_{\un} = L_{\un} U$ is generated by two conjugate regular unipotents. This in turn implies the result as in the $F_4(q)$-case.

Suppose that $q=9$. Then we consider another parabolic subgroup $P$ with $U$ of order $q^{31}$. In this case $L_{\un}\cong SL_3(q)\times  SL_2(q^2)$ \cite[Table 10]{C23}.
Then $L_{\un}$ is generated by two conjugate regular unipotents (see Remark \ref{SL29} and Lemma \ref{qq3}). 
The composition series of $L_{\un}$ on $U$ has four terms with factors of dimensions $3,4,12,12$ \cite[Table 2]{Ko14} which are factors of the lower central series of $U$. 
Arguing as above, we conclude that $P_\un$ is generated by two conjugate regular unipotents, and then the result follows.
\enp

The above results leave us with the cases where $G\in \{{}^2E_6(q),F_4(q)\}$, $q$ even. 
To consider the former case let $\mathbf{G}_{ad}$  be the simple algebraic group  of type $E_6$, $\mathbf{H}$  the simple algebraic group of type $F_4$, and we consider $\mathbf{H}$ as a subgroup of $\mathbf{G}$ defined
by $C_{\mathbf{G} } (\gamma)$, where $\gamma$ is a graph \au of order $2$. 

\bl{LS4}\cite[Theorem 5.1 and Remark following it]{LSe}
Let $\mathbf{G}$ be a simple algebraic group in characteristic p, $\si$ a Frobenius endomorphism of $\mathbf{G}$ and  $G=\mathbf{G}^\si$.   
\begin{itemize}
\item[$(1)$] All quasisimple subgroups $X$ of  $\mathbf{G}$ such that $X_{\un}/Z(X_{\un})\cong G_{\un}/Z(G_{\un})$   are conjugate  to $G$ in  ${\rm Aut}(\mathbf{G})$.  

\item[$(2)$] Let  $\tau=\si^k$, $k>1$ and $H= \mathbf{G}^\tau$. Then all quasisimple subgroups X of H such that $X_{\un}/Z(X_{\un})\cong G_{\un}/Z(G_{\un})$   are conjugate  to $G$ in ${\rm Aut}(H)$, in fact, in $N_{{\rm Aut}(\mathbf{G})}(H)$.
 
\item[$(3)$] {\rm \cite[Lemma 5.2]{LSe}} $N_{\mathbf{G}}(G_{\un})=G$.  
\end{itemize}\el

Recall that  the automorphism of $G_{\un}$ arising from the conjugation by elements of $N_{\mathbf{G}}(G_{\un})$ are called diagonal. 

Note that $(2)$ follows from $(1)$, see  \cite[Lemma 5.3]{LSe}. Our wording in the above statement  slightly differ from that in  \cite[Theorem 5.1]{LSe}.

If $\mathbf{G}$ is of adjoint type, then  $Z(G_{\un})=1$  and $N_{\mathbf{G}}(G_{\un}^\si)=G^\si$ \cite[Lemma 5.2]{LSe}. So in this case the statement of Lemma \ref{LS4} becomes simpler. In particular, this is the case for $\mathbf{G}$ of type $F_4$; in this case, $G=G_{\un}= N_{\mathbf{G}}(G_{\un})$. (Note that $\GG=F_4$ is of adjoint type.) 

\bl{Cu2}   
Let $G=F_4(q)$ and let $M$ be a maximal subgroup of $G$ containing $X\cong F_4(2)$. Then $M$ is conjugate to the standard subgroup of $G$ isomorphic to 
$F_4(q_1)$ with $q=q_1^r$, $r$ a  prime. Moreover, for every such $q_1$ there is exactly one subgroup containing X and isomorphic to $F_4(q_1)$.\el

\bp  The first claim is contained in  \cite{C23}. For the second one suppose the contrary, and let  $M,M_1$ be distinct subgroups containing $X$ and isomorphic to $F_4(q_1)$.
By Lemma \ref{LS4} or \cite{C23}, $gMg\up=M_1$ for some $g\in G$. Then $gXg\up <  M_1$.
Again by Lemma \ref{LS4},  $hgXg\up h\up =X$ for some $h\in M_1$. Then $hg\in N_G(X)<  N_{\GG}(X)$.  By Lemma \ref{LS4}(3), we have $N_{\GG}(X)=X$,
whence the claim. \enp

\bl{Cu3} The group $F_4(q)$, $q$ even, is generated by three conjugate regular unipotent elements.\el 

\bp By Lemma \ref{Cu2}, $F_4(q_i)$, $q_i^r = q=2^m$, where $r\mid m$ is a prime,
are the only maximal subgroups of $G=F_4(q)$ containing $X$. As $X$ is generated by
two conjugate regular unipotents $u,u'$, say, it suffices to show that
$G \setminus( \cup _i F_4(q_i))$ contains a conjugate  of $u$.
For this we show that the total number of conjugates of  $u$ in $G$ is
greater that the total number of conjugates of $u$ in $\cup _i F_4(q_i)$.
Note that the size the $G$-orbit of $u$ is $|G|/|C_G(u)|=|G|/q^4$ and
those in the above subgroups are $|F_4(q_i)|/{q_i}^4$.

We prove that $\frac{|G|}{q^4}>  \mathcal{S}(q)$, where
$\mathcal{S}(q)=\sum_{i} \frac{|F_4(q_i)|}{q_i^4}$ and the sum is taken
on all $2$-power $q_i$ such that $q=2^m = q_i^{r_i}$ with $r_i$ a prime.
Note that $2^{51m} <|F_4(q)|<2^{52m}$.
Hence, $$\mathcal{S}(q)\leq \sum_i q_i^{48}\leq
\sum_{\begin{smallmatrix} d\mid m\\ d\neq m\end{smallmatrix}}
(2^{d})^{48}\leq \sum_{i=0}^{\lfloor m/2\rfloor }
(2^i)^{48}=\frac{(2^{48})^{1+\lfloor m/2\rfloor} -1}{2^{48}-1} < 2^{25m}
<\frac{|G|}{q^4}.$$
\enp

Let $\mathbf{G}$ be a simple algebraic group of type $E_6$, and let $\si_j$,   $j=1,2$, be a Frobenius morphism $\mathbf{G}\ra \mathbf{G}$ such that $G_k:=\mathbf{G}^{\si^k_1}\cong E_6(2^k)$,  ${}^2G_k:=\mathbf{G}^{\si^k_2}\cong {}^2E_6(2^k)$. We choose $\si_k^1$ to be a standard Frobenius arising from $x\mapsto x^{2^k}$ for $x\in\overline{\FF}_2$,
and $\si_2^k=\gamma\cdot \si^1_k=\si_1^k\cdot \gamma$. (Note that $\gamma$ commutes with $\si_1^k$ and $\si_2^k$.) We call $G_k$ and ${}^2G_k$ standard finite subgroups of type $E_6$ and ${}^2E_6$, respectively,  of~$\mathbf{G}$.

Set $H^+_k=C_{\mathbf{G}}(\gamma, \si_1^k)$ and $H^-_k=C_{\mathbf{G}}(\gamma, \si_2^k)$. We can apply $\gamma$ to $\GG$  and then $\si_j^k$ to $C_{\GG}(\gamma)$ or conversely, first $\si_j^k$ to $ \GG $ and then $\gamma$ to $C_{\GG}(\si_j^k)$, to obtain $H^+_k$ and $H^-_k$. 
Therefore,    $H^+_k$ is the intersection of $ C_{ \mathbf{G}}(\si_1^k)=E_6(2^k)$ 
and $\mathbf{H}:= C_{\GG}(\gamma) $ and  $H^-_k$ is the intersection of 
$ C_{ \mathbf{G}}(\si_2^k)={}^2E_6(q^k)$ and $\mathbf{H}$. 
We call $H^+_k$ and $H^-_k$ standard finite subgroups of type $F_4$ in $ \mathbf{G}$. 

In Lemmas \ref{Cum}, \ref{q66} and  \ref{3co} we deals with groups $G={\mathbf G}_{ad}^\si$ of adjoint type. In this case $G$ is not always simple whereas $G_{\un}$ is a simple group of Lie type. For our purposes it is convenient to call a subgroup $M$ of $G$ maximal if $M\cap G_{\un}$ is a proper maximal subgroup of $G_{\un}$. The conjugacy classes of  maximal subgroups of $G$ and $G_{\un}$ are determine in \cite{C23}.

\bl{Cum}\cite{C23}
Let $G={}^2E_6(q)_{ad}$ and let $M$ be a maximal subgroup of $G$ containing $X\cong F_4(2)$. Then $M$ is conjugate to the standard subgroup of $G$ of type $H_m\cong  F_4(q)$ or of type ${}^2E_6(q_1)_{ad}$ with $q=q_1^r$, $r$ an odd prime. 
\el

\bp 
Note that $X <  G_{\un}$. By inspection of the list of  maximal subgroups of $G$ in  \cite{C23}, one observes that those   containing a subgroup isomorphic to $F_4(2)$ 
are as indicated in the lemma. In addition, groups $M$ isomorphic to each other  are conjugate in $G$, that is, form a single conjugacy class.
So one can  choose the standard subgroup as a representative of each class.\enp

\bl{q66}
Let $X$ be a subgroup of $G={}^2E_6(q)_{ad}$ isomorphic to $H_k^- $. Then $X$ is conjugate in $G$ to  $H_k^- $. In addition, $N_{G}(H_k^-)=H_k^-$.\el

\bp   Let $q=2^m$,   so $G={}^2G_m$. Then $X<N_G(X)\leq M$,   
where $M$ is a maximal subgroup of $G$.  As $G$ is of adjoint type, by Lemma \ref{Cum}, we can assume that $M$ 
is the standard subgroup of $G$, so       $M=H^-_m$ or   $M={}^2G_{m_1}$,  
where $q=2^{m_1r}$ with  $r$ an odd prime. If  $M=H^-_m$  then, by Lemma \ref{LS4}, 
$X$ is conjugate in $M$ to a standard subgroup of $H^-_m$. As standard subgroups of $H^-_m$ are standard in $G$, the result follows. In addition, this implies the result if $m$ is a $2$-power as in this case $G$ has no subgroup of type ${}^2 E_6(q_1)$ for $q_1<q$.  

Suppose that  $M={}^2G_{m_1}$. Then the result follows by induction   on $m$; the base of induction is $m_2$, the   $2$-power part of $m$, case already settled. 

The additional case follows from Lemma \ref{LS4}(3) if $N_G(X)\leq M=H_m^-$, otherwise this follows by induction.\enp

\bl{3co} Let $X=F_4(2)<M<G={}^2E_{6}(q)_{ad}$, $q$ even, where  $M$  is a maximal subgroup of $G$ not containing $G'$.
Suppose that  $X < gMg\up$ for $g\in G$. Then $g\in  M$. Consequently, if 
$M,M'$ are isomorphic maximal subgroups of $G$ containing $X$ then $M=M'$. 
\el

\bp 
Let $q=2^m$. By Lemma \ref{Cum},  
$M$ is conjugate to a standard subgroup  $H_m\cong F_4(q)$ or ${}^2G_{m_1}\cong {}^2G_{ad}(q_1)$ with $q=q_1^r$, where $r=m/m_1$ is  an odd prime. 
We can assume  that $M$ itself is standard.

Suppose the contrary, that $X<gMg\up\neq M$. Then  $g\up Xg<M$. By Lemma \ref{q66}, $X$ and $g\up Xg $ are conjugate in $M$ so $X=y\up g\up Xgy$ for some $y\in M$. So $gy\in N_G(X)$.  As $N_G(X)=X$ by Lemma \ref{q66}, 
we have $gy\in M$ and $g\in M$, a contradiction. 

The complementary assertion follows from this as isomorphic maximal subgroups of $G$ (containing $X$) are conjugate by \cite{C23}. \enp

\bl{ng7} Let $G={}^2E_6(q)_{ad}$, $q$ even. Distinct maximal subgroups of $G_{\un}$ containing a fixed subgroup $X\cong F_4(2)$ are not isomorphic.
Consequently, if $q=2^m$ then  the number of maximal subgroups of $G_{\un}$ containing $X$ is $d+1$, where $d$ is the number of  odd primes, dividing $m$.  \el

\bp 
We only need to consider the case with $G_{\un}\neq G $ by Lemma~\ref{3co}. This   implies that $3\mid (q+1)$  and $|G/G_{\un}|=3$.  

Let $M$  be a maximal subgroup of $G$, and $X<M.$  
If $M\cong F_4(q)$ then the reasoning in Lemma~\ref{3co} works (in this case $M=M_{ad}=M_{\un}$ so the subgroups of $M$ isomorphic to $X$ are conjugate in~$M$).

Let $M={}^2E_6(q_1)_{ad}<G$, $q=q_1^3$. Then $3 \mid (q_1+1)$ and 
 $X<M_{\un}<M$. As $M_{\un}$ is normal in $M$,  $M_{\un}$ is the unique subgroup of $G_{\un}$ containing $X$. 

The assertion on the number of maximal subgroups in question follows from \cite{C23}. 
\enp

\bl{2e6q}
Let $G={}^2E_6(q)_{ad}$, $q$ even. Then  $G_{\un}$ is generated by three conjugate (in $G_{\un}$) regular unipotent elements.
\el

\bp  
Let $X<G_{\un}$, $X\cong F_4(2)$ be as above. Let $q=2^m$, and let $s$ be 
the number of distinct odd prime divisors of $m$.

By Lemma \ref{ng7}, $F_4(q)$ and ${}^2E_6(q_i)_{\un}$, $q_i^r=q$, where $r\mid m$
is an odd prime, are the only maximal subgroups generated by unipotent elements and
containing $X$. As $X$ is generated by two conjugate regular unipotents $u,u'$, say, it suffices to show that  
${}^2E_6(q)_{\un} \setminus(F_4(q)\cup _i {}^2E_6(q_i)_{\un})$, contains a conjugate  of $u$. For this we show that the total number of conjugates of  $u$ in $G_\un$ is greater that the total number of conjugates of $u$ in  
$F_4(q)\cup _i {}^2E_6(q_i)_{\un}$.

Note that the size the $G_\un$-orbit of $u$ is $|G_\un|/|C_{G_\un}(u)|=|G_\un|/q^6$,
and those in the above subgroups are are $|F_4(q)|/q^4$ and $|{}^2E_6(q_i)_\un|/q_i^6$.

We prove   that $\frac{|G_\un|}{q^6}>\frac{|F_4(q)|}{q^4}+\sum_i\frac{|{}^2E_6(q_i)_\un|}{q_i^6}$. 
Set $\mathcal{S}(q)=\sum_{i} \frac{|{}^2E_6(q_i)_{\un}|}{q_i^6}$, 
where the sum is taken on all $2$-power $q_i$ such that $q=2^m=q_i^{r_i}$ with $r_i$ odd prime.  

Assume $m>1$. First of all, note that $2^{77m}<|E_6(q)_{\un}|<2^{78m}$ and $|F_4(q)|< 2^{52m}$.
Now, $$\mathcal{S}(q)\leq \sum_i q_i^{72}\leq \sum_{\begin{smallmatrix} d\mid m\\ d\neq 1\end{smallmatrix}} (2^{m/d})^{72}
= \sum_{\begin{smallmatrix} d\mid m\\ d\neq m\end{smallmatrix}} (2^{d})^{72}\leq \sum_{i=0}^{\lfloor m/2\rfloor } (2^i)^{72}=\frac{(2^{72})^{1+\lfloor m/2\rfloor} -1}{2^{72}-1} \leq 2^{37m}.$$ Hence, 
$$ \frac{|F_4(q)|}{q^4}+\mathcal{S}(q)<2^{48m}+2^{37m}< 2^{71m}<\frac{|{}^2E_6(q)_{\un}|}{q^6}. $$

If $m=1$ (or $m$ is a $2$-power) then we only need to observe that $|G_{\un}|>q^2|F_4(q)|$,
which is clear from the above.  \enp

Note that the method used for $G={}^2E_6(q)$ works for $E_6(q)$. However, we used for this group an alternative approach. 

In conclusion,  Theorem \ref{t33} follows from the previous results according to the Table \ref{casi}.

\begin{table}[htp]
\begin{center} 
\begin{tabular}{|ll|}\hline
$SL_2(q)$, $q\geq 4$ even  & \cite[Lemma 3.1]{GS} \\
$SL_2(9)$ & \cite[Lemma 3.1]{GS} \\
$SL_n(q)$, $n\geq 3$ &  Lemma \ref{+23}\\\hline
$SU_4(q)$, $q\geq 4$ even & Lemma \ref{uu4}\\
$SU_5(q)$, $q\geq 5$ odd  &  Lemma \ref{unn}\\
$SU_5(q)$, $q\geq 4$ even  &  Lemma \ref{uu5}\\
$SU_n(q)$, $n\geq 7$ odd &  Lemma \ref{unn}\\\hline
$Sp_4(q)$, $q$ even & Lemma \ref{s44}\\
$Sp_{2n}(q)$, $n\geq 3$ & Lemma \ref{spp2} \\\hline
$\Omega_{2n+1}(q)$, $n\geq 3$ and $q$ odd & Lemma \ref{oo2}\\
$\Omega_{2n}^+(q)$, $n\geq 4$ and $q$ even & Lemma \ref{oo1}\\
$\Omega_{2n}^-(q)$, $n\geq 4$ & Lemma \ref{oo2} \\\hline
$G_2(q)$, $q\geq 4$ even & Lemma \ref{G2even}\\
$F_4(q)$  & Lemmas \ref{f44}  and \ref{Cu3} \\ \hline
${}^2E_6(q)$ & Lemmas  \ref{2e6} and \ref{2e6q}\\ \hline
\end{tabular}
\end{center}
\caption{Auxiliary results for Theorem \ref{t33}.}\label{casi}
\end{table}

Theorem \ref{p33} follows from Lemmas \ref{+23}, \ref{oo1}, \ref{oo2}, \ref{pa2}, \ref{ee6} and \ref{2e6}.
Finally, Theorem \ref{qs3} follows from Theorem \ref{t22}, Theorem  \ref{t33} and  Lemmas \ref{unn}, \ref{oo1}, 
and \ref{ee6}.

\section*{Acknowledgments}

The first author is partially supported by INdAM-GNSAGA.

\end{document}